\documentclass[11pt]{article}

\usepackage{graphicx}
\usepackage{latexsym,amsmath,amsfonts,amscd, amsthm, dsfont}
\usepackage{color}
\usepackage{epsfig,verbatim,epstopdf,graphics}
\usepackage{changebar}
\usepackage{scalerel}
\usepackage{chngcntr}

\counterwithin*{equation}{section}
\graphicspath{{./}{./figure/}}

\newtheoremstyle{plainNoItalics}{}{}{\normalfont}{}{\bfseries}{.}{ }{}

\theoremstyle{plain}
\newtheorem{thm}{Theorem}[section]

\theoremstyle{plainNoItalics}

\newtheorem{rem}[thm]{Remark}

\newcommand{\bE}{\mathbf{E}}
\newcommand{\bB}{\mathbf{B}}
\newcommand{\bJ}{\mathbf{J}}
\newcommand{\bA}{\mathbf{A}}
\newcommand{\bn}{\mathbf{n}}
\newcommand{\bw}{\mathbf{w}}
\newcommand{\bx}{\mathbf{x}}
\newcommand{\by}{\mathbf{y}}

\newcommand{\bT}{\mathbf{T}}
\newcommand{\bR}{\mathbf{R}}

\newcommand{\bj}{\mathbf{j}}
\newcommand{\bk}{\mathbf{k}}
\newcommand{\bl}{\mathbf{l}}
\newcommand{\bg}{\mathbf{g}}
\newcommand{\bgamma}{\mathbf{\gamma}}

\newcommand{\be}{\begin{eqnarray}}
\newcommand{\ee}{\end{eqnarray}}
\newcommand{\beno}{\begin{eqnarray*}}
\newcommand{\eeno}{\end{eqnarray*}}

\makeatletter

\newcommand{\Rmnum}[1]{\expandafter\@slowromancap\romannumeral #1@}
\makeatother

\begin{document}

\baselineskip=1.8pc

%\vspace*{.10in}

%=============  title  =========================
\begin{center}
{\bf
An Asymptotic Preserving Maxwell Solver Resulting in the Darwin Limit of Electrodynamics
%Coupled with Semi-Lagrangian Finite Difference WENO for Vlasov Simulations
%Kinetic Simulations and Advection in Incompressible Flows
}
\end{center}

\vspace{.2in}
\centerline{
Yingda Cheng\footnote{
Department of Mathematics, Michigan State University, East Lansing, MI, 48824. E-mail: ycheng@math.msu.edu. Research is supported by NSF grant  DMS-1318186.
},
Andrew J. Christlieb\footnote{
Department of Computational Mathematics, Science and Engineering, and Department of Mathematics, Michigan State University, East Lansing, MI, 48824. E-mail: andrewch@math.msu.edu
},
 Wei Guo\footnote{
Department of Mathematics, Michigan State University, East Lansing, MI, 48824. E-mail:
wguo@math.msu.edu},
Benjamin Ong\footnote{Department of Mathematical Sciences, Michigan Technological University, Houghton, MI, 49931. E-mail:
 ongbw@mtu.edu.}
}

\bigskip
\noindent
\centerline{{\bf Abstract}}
%\bigskip

 In plasma simulations, where the speed of light divided by a  characteristic length  is at a much higher frequency than other relevant parameters in the underlying system, such as the plasma frequency, implicit methods begin to play an important role in generating efficient solutions in these multi-scale problems.  Under conditions of scale separation, one can rescale Maxwell's equations in such a way as  to give a  magneto static limit known as the Darwin approximation of  electromagnetics.
%
  %This model converts Maxwell's equations into a set of elliptic equations.
 In this work, we present a new approach to solve Maxwell's equations based on a Method of Lines Transpose (MOL$^T$) formulation, combined with a fast summation method with computational complexity $O(N\log{N})$, where $N$ is the number of grid points (particles). Under appropriate scaling, we  show that the proposed schemes result in   asymptotic preserving methods that can recover the Darwin limit of electrodynamics.
\vfill

{\bf Keywords:} implicit method; Maxwell's equations; Darwin approximation; Method of Lines Transpose; fast summation method; asymptotic preserving method.
%\textcolor{red}{please do a thorough spell check}

\newpage

\section{Introduction}
\label{sec:introduction}

In this paper, we develop asymptotic preserving (AP) numerical methods for   Maxwell's equations which can recover the Darwin limit of electrodynamics under an appropriate scaling limit.
AP schemes are also known as temporal multi-scale methods in the literature, see \cite{jin2010asymptotic} for review and recent development of such methods. The main attractive features lie in the ability to preserve the asymptotic limit of the underlying equation at the discrete level, permitting large time step evolution even when the scaling parameter becomes small in the system. There has been lots of success of AP schemes in the kinetic simulations, in which the Knudsen number $\varepsilon$ is used to characterize the kinetic scales. As  $\varepsilon$ goes to 0, the AP scheme is  designed in a way that it becomes a consistent discretization of the limiting hydrodynamic models, which, in fact, is mimicking the asymptotic limiting procedure on the partial differential equation (PDE) level. As a result, the AP scheme is uniformly stable with respect to $\varepsilon$. Recent development of AP schemes  for kinetic models include a penalty method   for solving the Boltzmann equation with general collisional operator \cite{filbet2010class}, the macro-micro decomposition technique for  the BGK model  \cite{bennoune2008uniformly}, among many others \cite{pieraccini2007implicit,pieraccini2012microscopically,xiong2015high}. Based on similar  ideas, AP schemes were also developed for the Euler-Poisson and Euler-Maxwell models in \cite{degond2008analysis,degond2012numerical} to recover the quasi-neutral limit of the incompressible Euler model as   the Debye length $\lambda_D$ goes to 0.

%in which, under a suitable scaling, the Debye length $\lambda_D$ is incorporated to characterize the quasi-neutrality of the plasma system. Such a method is able to capture the quasi-neutral limit as $\lambda_D$ goes to 0, in the sense that the scheme becomes a consistent discretization for the incompressible Euler model. Similarly to the kinetic model, the AP property of the schemes ensures the uniform stability with respect to $\lambda_D$, which allows time steps much larger than the electron time scale.

This paper, on the other hand,  concerns the construction of numerical schemes for Maxwell's equations which can capture asymptotic limit to the Darwin model. The Darwin model is a well-known approximate model to Maxwell's equations \cite{darwin1920li}. It is obtained from Maxwell's equations  by neglecting the solenoidal, i.e., the divergence-free part of the displacement current in Amp\`{e}re's law. This results in a set of  elliptic equations with separated electric and magnetic fields which are easier to solve than full Maxwell's equations. In \cite{degond1992analysis}, Degond \emph{et al.} showed that the Darwin model approximates Maxwell's equations up to second order for the magnetic field and   third order for the electric field with respect to the dimensionless parameter $\epsilon=\frac{\bar{v}}{c}$ in a three-dimensional bounded simply connected domain, where $\bar{v}$ is a characteristic velocity and $c$ is the speed of light. Such an analysis in fact verifies the effectiveness of the Darwin model when no high frequency phenomenon or rapid change occurs in the physics system. Many research efforts have been devoted to the development of numerical schemes for solving the Darwin model. For instance,
 in \cite{ciarlet1997finite}, a finite element method was proposed and the well-posedness of the associated variational problems was also established. In plasma physics, numerical schemes to the Vlasov-Darwin model, which is a simplified model for the Vlasov-Maxwell system, have been considered in \cite{sonnendrucker1995finite, schmitz2006darwin, besse2007numerical}. In \cite{raviart1995approximate, raviart1996hierarchy}, a hierarchy of approximate models for Maxwell's equations are established with the perfectly conducting and the Silver-M\"{u}ller boundary conditions. The quasistatic Darwin models are proven to be first and second order approximation to Maxwell's equations with respect to $\epsilon$. It is therefore, of great interest to applications, to construct efficient Maxwell solvers that can capture the Darwin limit automatically.

%Our focus in the paper is to study the asymptotic behavior of several potential formulations of Maxwell's equation, including a vector potential which was firstly introduced in \cite{de1998time} and the standard scalar and vector potentials with the Coulomb and Lorentz gauge.

The approach we use in this paper consists of several important components, including the development and the extension of the  Method of Line Transpose (MOL$^T$) framework, the investigation of the scalar and vector potential formulations of Maxwell's equations and their asymptotic limit as key steps  to ensure the solver to capture the correct Darwin limit,  and an efficient treecode algorithm to further accelerate the computation.  The MOL$^T$ method we consider in this paper is also known as transverse MOL, and Rothe's method in the literature \cite{salazar2000theoretical,schemann1998adaptive}. As the name implies, discretization is carried out in an orthogonal fashion, where the time variable is first discretized, followed by solving the resulting boundary value problems (BVPs) at discrete time levels. %However, the MOL$^T$ approach was mostly overlooked in history, partially due to the extra complications for efficiently solving BVPs, especially when the boundary layers appear in the solution. On the other hand, under the integral method framework,
The MOL$^T$ approach is advantageous when coupled with the integral method framework since one can employ many fast summation methods, such as the fast multipole method  \cite{greengard1987fast} (FMM) and the treecode \cite{barnes1986hierarchical} to reduce the computational complexity of evaluation from $O(N^2)$ to $O(N)$ or $O(N\log{N})$. In \cite{causley}, the MOL$^T$ method was developed for the wave equations, in which the BVP is first split into a series of one-dimensional BVPs via the alternating direction implicit (ADI) technique, then the proposed one-dimensional solver is applied to the split BVPs at a price of splitting errors. The resulting method is A-stable, easy to implement and  the computation complexity can be reduced to $O(N)$ by utilizing the analytical properties of the one-dimensional Green's function, see \cite{causley2013method}. Arbitrary temporal accuracy is attained by successive convolutions in \cite{causley2014higher}. However, there are still several challenges  for the extension of this method to Maxwell's equations. To employ the framework for the wave equation, first we need the potential formulations of Maxwell's equations. They also turn out to be crucial to achieve AP properties for the schemes. Second, the ADI splitting strategy used in \cite{causley} is no longer suitable when the Silver-M\"{u}ller boundary condition is imposed.   Instead, we  shall invert the three-dimensional Helmholtz operator. Note that there is a huge amount of literature on integral methods for Maxwell's equations in the frequency domain, where very fast algorithms have been developed, such as the FMM \cite{darve2000fast,otani2008periodic,nishimura2002fast}. However, the FMM deals with the oscillating free-space Green's function.
%	 in the form of $\exp(-ikr)/r$ where $r=|\mathbf{x}-\mathbf{y}|$ and $k$ is a real number.
	   While, in our case,  the associated Green's function exhibits exponential decay and has been recently incorporated in both FMM \cite{gimbutas2003generalized} and treecode \cite{li2009cartesian} algorithms. We use the treecode algorithm to speed up our calculation, and  similar to the scheme in \cite{causley}, the newly proposed scheme is  unconditionally stable due to the implicit treatment in the MOL$^T$ framework, which is a highly desirable property in the plasma simulations, since $\frac{1}{\epsilon}$ is at a much higher frequency than the relevant parameters in the underlying system, such as the plasma frequency. In \cite{mavsek2010mesh,gibbon2010progress}, the treecode is used to solve the Darwin model in plasma simulations.

 The rest of the paper is organized as follows. In Section \ref{sec:models}, we introduce the underlying models including Maxwell's equations and the Darwin model. In particular, we consider Maxwell's equations which are written in terms of  potentials. We show that, under a suitable scaling, the potential forms of Maxwell's equations are consistent with the Darwin model up to certain orders of the dimensionless parameter $\epsilon$. In Section \ref{sec:semidiscrete}, we formulate the semi-discrete schemes for solving the rescaled potential forms of Maxwell's equations in the MOL$^T$ framework. Through   formal asymptotic analysis, we  show that the   semi-discrete schemes are AP in the sense that the schemes can capture the Darwin limit as $\epsilon$ goes to 0. In Section \ref{sec:fullydiscrete}, we propose several fully discrete schemes and utilize the treecode to speed up the computations. Two numerical examples are presented in Section \ref{sec:numerical} to verify the performance of our methods. We end with concluding remarks  and future work in Section \ref{sec:conclusion}.

\section{The Models}
\label{sec:models}

In this section, we review Maxwell's equations, their potential formulations, and the asymptotic limit to the Darwin model.

\subsection{Maxwell's Equations and the Potential Formulations}
\label{sec:equation}
We are interested in Maxwell's equations defined on $\Omega \subset \mathds{R}^3$, which can be written in MKS units as follows:
\begin{subequations}
\label{eq:maxwell}
\begin{align}
&-\frac{1}{c^2} \frac{\partial  \bE}{\partial t} +\nabla \times \bB =  \mu_0 \bJ,\label{eq:maxwell_e}\\
&\frac{\partial  \bB}{\partial t} +\nabla \times \bE =  0,\\
&\nabla \cdot \bE =\frac{\rho}{ \epsilon_o },  \\
& \nabla \cdot \bB = 0,
\end{align}
\end{subequations}
subject to the continuity equation
\be
\label{continuity}
\frac{\partial \rho}{\partial t}+\nabla \cdot \bJ=0,
\ee
where $c$ is the speed of light,   $\mu_0$ and $\epsilon_o$ are the magnetic permeability and the electric permittivity, respectively. Note that $\mu_0$ and $\epsilon_o$ are related to $c$ according to $\mu_0\epsilon_oc^2=1$.

We consider two types of boundary conditions: the perfectly conducting boundary condition
\begin{equation}
\bE \times \bn=0 \quad \textrm{on} \, \,\Gamma_C,
\label{eq:pec_bc}
\end{equation}
where $\bn$ is outward unit normal vector on the boundary,
and the Silver-M\"{u}ller boundary conditions
\begin{equation}
\label{eq:sl_bc}
  (\bE-c\bB\times \bn)\times \bn=\bg \times \bn  \quad \textrm{on} \, \,\Gamma_A,
\end{equation}
where $\bg$ is a given function. Note that in the case of $\bg=0$, the Silver-M\"{u}ller boundary conditions correspond to the absorbing boundary condition.
Here, $\Gamma_C$ and $\Gamma_A$ denote the subsets of the boundary on which the perfect conducting and Silver-M\"{u}ller boundary conditions are imposed.

Next, we review several potential formulations of Maxwell's equations for the purpose of developing our numerical schemes. First, we consider the potential formulation    proposed in \cite{de1998time}. Let $\bw$ be the time integral of the electric field, i.e., $\bw(t, \bx)= \int_0^t \bE(s, \bx) ds$. A useful property of $\bw$ is that  both the electric and magnetic fields can be represented by such a single vector potential. Therefore, it is relatively easy to impose boundary conditions. In particular, substituting $\bw(t, \bx)= \int_0^t \bE(s, \bx) ds$ into \eqref{eq:maxwell}, we obtain a set of evolution equations for $\bw$ as follows.
\begin{subequations}
\label{new}
\begin{align}
&\frac{1}{c^2} \frac{\partial^2  \bw}{\partial t^2} +\nabla \times (\nabla \times \bw) =  \nabla \times \bB(0, \bx)-\mu_0 \bJ,\\
&\nabla \cdot \left (\frac{\partial  \bw}{\partial t} \right)=\frac{\rho}{\epsilon_0},\\
& \bw(0, \bx)=0,\\
&\partial_t \bw(0, \bx)=\bE(0, \bx).
\end{align}
\end{subequations}
Note that
the electric and magnetic fields can be conveniently expressed in terms of $\bw$ as $\bE(t, \bx)=\partial_t \bw(t, \bx)$, and $\bB(t, \bx)=-\nabla \times \bw(t, \bx)+\bB(0, \bx)$. In other words, one can just solve for $\bw$ and obtain $\bE$ and $\bB$ by a numerical differentiation procedure. %Moreover, recall that $\bw$ is simply the time integral of the electric field, and hence the boundary conditions for $\bw$ could be naturally imposed. In particular,
The perfectly conducting boundary condition \eqref{eq:pec_bc} becomes
\begin{equation}
\label{eq:pec_bc_w}
\bw \times \bn=0 \quad \textrm{on} \, \,\Gamma_C,
\end{equation}
and the Silver-M\"{u}ller boundary conditions \eqref{eq:sl_bc} become
\begin{equation}
\label{eq:sl_bc_w} \left (\partial_t \bw + c \,(\nabla \times \bw) \times \bn \right ) \times \bn=\left (c \, \bB(0,\bx)
\times \bn+\bg \right) \times \bn   \quad \textrm{on} \, \,\Gamma_A.
\end{equation}

%\textcolor{red}{Another definition below?}
%\emph{Since $\nabla \cdot \bB(0,x)=0$, there exists $\bA$, such that $\bB(0,x)=\nabla \times \bA$. Then define $\bw(t, x)= \int_0^t \bE(s, x) ds-\bA$, we have $\bE(t, x)=\partial_t \bw(t, x)$, and $\bB(t, x)=-\nabla \times \bw(t, x)$, and $$\frac{1}{c^2} \frac{\partial^2  \bw}{\partial t^2} +\nabla \times (\nabla \times \bw) =  -\mu_0 \bJ$$ and we have get rid of the term relating to $\bB(0,x)$.}

We remark that equations \eqref{new} also have an equivalent wave formulation:
\begin{subequations}
\label{new2}
\begin{align}
&\frac{1}{c^2} \frac{\partial^2  \bw}{\partial t^2} -\Delta \bw=
-\nabla\left(\int_0^t \frac{\rho}{\epsilon_0} ds\right) +\nabla \times \bB(0, \bx)-\mu_0 \bJ,\\
&\nabla \cdot \left (\frac{\partial  \bw}{\partial t} \right)=\frac{\rho}{\epsilon_0},\\
& \bw(0, \bx)=0,\\
&\partial_t \bw(0, \bx)=\bE(0, \bx).
\end{align}
\end{subequations}
%Note that each component of $\bw$ satisfies a decoupled wave equation, but
Notice that here the source term involves the time integral of the density function, and the constraint $ \nabla \cdot \left (\frac{\partial  \bw}{\partial t} \right)=\frac{\rho}{\epsilon_0}$ is necessary to enforce the continuity equation. The numerical scheme discussed in this paper  will be based on the wave equation formulation \eqref{new2}. A numerical scheme for \eqref{new}  will be considered in future work.%We will also comment on the curl-curl formulation \eqref{new} and leave the implementation of scheme to  future study.

Besides the vector potential formulation introduced above, we can also consider the more common scalar $\psi$ and vector potential $\bA$ defined by
$$\bE=-\nabla \psi-\frac{\partial \bA}{\partial t}, \quad \bB=\nabla \times \bA.$$
Here,  gauge conditions are important to uniquely determine the potentials $\psi$ and $\bA$. %In fact, enforcing the gauge conditions are equivalent to enforcing the continuity condition.
For example, with the Lorentz gauge
\be
\label{lorentzgauge}
\nabla \cdot \bA+\frac{1}{c^2}\frac{\partial \psi}{\partial t}=0,
\ee
Maxwell's equations can be written as a set of wave equations for the potential $\psi$ and $\bA$:
\begin{subequations}
\label{lorentz}
\begin{align}
&\frac{1}{c^2} \frac{\partial^2  \psi}{\partial t^2} -\Delta  \psi =  \frac{\rho}{\epsilon_0},\\
&\frac{1}{c^2} \frac{\partial^2  \bA}{\partial t^2} -\Delta  \bA =  \mu_0 \bJ.
\end{align}
\end{subequations}
This formulation is especially suited for a particle code, where $\rho$ and $\bJ$ are linear combinations of  Dirac delta functions.
Another common choice is to take the Coulomb gauge
\be
\label{coulombgauge}
\nabla \cdot \bA=0.
\ee
Under this gauge, the evolution equations for the potentials can be written as
\begin{subequations}
\label{coulomb}
\begin{align}
&  -\Delta \psi =  \frac{\rho}{\epsilon_0},\label{eq: coulomb1}\\
&\frac{1}{c^2} \frac{\partial^2  \bA}{\partial t^2} -\Delta  \bA =  \mu_0 \bJ-\frac{1}{c^2}\nabla \left(\frac{\partial \psi}{\partial t}\right),\label{eq: coulomb2}
\end{align}
\end{subequations}
where the first equation \eqref{eq: coulomb1} is simply Poisson's equation for the scalar potential $\psi$, and  the second equation \eqref{eq: coulomb2} is a set of wave equations for each component of the vector potential $\bA$. Note that \eqref{eq: coulomb2} involves $\nabla \left(\frac{\partial \psi}{\partial t}\right)$, which is not simple to compute for a particle code.

%However, the coupled boundary conditions can be obtained as follows.
As for the boundaries,  we have
\begin{equation}\label{eq:bc_per_A}
  (\nabla\times\bA) \cdot\bn = 0,
\end{equation}
from $\bB\cdot\bn=0$. If $\bA$ is obtained, the perfectly conducting boundary condition becomes
\begin{equation}\label{eq:bc_per_psi}
  \nabla \psi\times\bn= -\partial_t \bA\times\bn \quad \textrm{on} \, \,\Gamma_C,
\end{equation}
 from $\bE\times\bn=0$.
 For the Silver-M\"{u}ller conditions, by replacing $\bE$ and $\bB$ in \eqref{eq:sl_bc} with $\psi$ and $\bA$, we have
 \begin{equation}
 (\nabla\psi+\partial_t\bA+c\nabla\times\bA\times \bn)\times \bn=-\bg \times \bn \quad \textrm{on} \, \,\Gamma_A.
 \end{equation}

 Comparing with the potential formulation for $\bw$, we can see that   it is non-trivial to impose   boundary conditions on $\psi$ and $\bA$ in a \emph{decoupled} manner.
Obviously,  the  coupling between $\psi$ and $\bA$ poses implementation challenges. Therefore, in this paper, we only implement the scheme derived from the $\bw$ formulation, and leave the implementation of the $\psi, \bA$ formulation to future work. However,   for completeness, we   provide theoretical analysis for both potential formulations in Sections \ref{sec:models} and \ref{sec:semidiscrete}.

\subsection{The Darwin Model}
 In this subsection, we review the Darwin model while omitting the details of derivation. The readers are referred to \cite{degond1992analysis} for more discussions about the model  and its well-posedness.
Consider the Helmholtz decomposition, $\bE=\bE_{irr}+\bE_{sol}$, where $\nabla \times \bE_{irr}=0$ and $\nabla \cdot \bE_{sol}=0,$ i.e.,
  $\bE_{irr}$ and $\bE_{sol}$ refer to the irrotational and solenoidal components of
 the electric field $\bE$, respectively. Then, we have the following  elliptic equations.
 \beno
&&\-\Delta \psi =\frac{\rho}{\epsilon_0}, \quad \bE_{irr}=-\nabla \psi, \\
&& \psi|_{\Gamma_{C_i}}=\alpha_i, \quad 0 \le i \le m,
\eeno
where $\Gamma_{C_i}$ denote the connected components of the boundary $\Gamma_C$, and $\Gamma_{C_0}$ denotes the outer boundary. $\{\alpha_i, \, i=0, \ldots m\}$ are solutions of the differential system
 \beno
&&\sum_{j=0}^m c_{ij} \frac{d \alpha_j}{dt}=\frac{1}{\epsilon_0} \int_\Omega \bJ \cdot \nabla \chi_i \, dx, \\
&& \alpha_i(t=0)=\alpha_{i0}, \quad 0 \le i \le m.
\eeno
Here $\{\chi_i, \, i=0, \ldots m\}$ are the solutions of
 \beno
&&\Delta \chi_i=0, \\
&& \chi_i|_{\Gamma_{C_i}}=\delta_{ij}, \quad 0 \le j \le m,
\eeno
and $\alpha_{i0}$ depends on the initial value of $\bE_{irr}$, and the capacitance coefficients $c_{ij}=\int_{\Gamma_{C_j}} \frac{\partial \chi_i}{\partial \bn} ds.$

The Darwin model can be derived from Maxwell's equations by neglecting the transverse component $\frac{\partial\bE_{sol}}{\partial t}$ of $\frac{\partial\bE}{\partial t}$ in \eqref{eq:maxwell_e}. In particular,
the magnetic field $\bB$ satisfies
\beno
&& \nabla \times \nabla \times \bB=\mu_0 \nabla \times \bJ, \quad \nabla \cdot \bB=0, \\
&& \bB \cdot \bn= \bB(0,\bx) \cdot \bn, \quad (\nabla \times \bB) \times \bn=\mu_0 \bJ \times \bn \quad \textrm{on} \, \, \Gamma_C.
\eeno
$\bE_{sol}$ satisfies
\beno
&& \nabla \times \nabla \times \bE_{sol}=-  \frac{\partial}{\partial t} \nabla \times \bB, \quad \nabla \cdot \bE_{sol}=0.\\
&& \bE_{sol} \times \bn=0\quad\textrm{on} \, \, \Gamma_C, \quad \int_{\Gamma_{C_i}}\bE_{sol} \cdot \bn ds=0  \quad 0\leq i\leq m.
\eeno

Imposing the  Silver-M\"{u}ller boundary conditions is more complicated, which was discussed in \cite{raviart1996hierarchy}. For example, $\bB$ satisfies
$$ (c^2 \nabla \times \bB- c \partial_t\bB \times \bn)\times \bn= \frac{1}{\epsilon_0} \bJ \times \bn +\partial_t \bg \times \bn \quad \textrm{on} \, \, \Gamma_A,$$
while $\bE_{sol}$ satisfies
$$ (\partial_t \bE_{sol}+c (\nabla \times \bE_{sol}) \times \bn)\times \bn= \partial_t \bg \times \bn \quad \textrm{on} \, \, \Gamma_A.$$

\subsection{Asymptotic Analysis of the Models}

The goal of this paper is to develop  implicit AP Maxwell solvers that can recover the Darwin limit under appropriate scaling limit. To serve the purpose, we will review the asymptotic analysis of Maxwell's equations and the potential formulations, and establish their connection with the Darwin model in this subsection.

\subsubsection{A Scaling}
\label{sec:scaling}
First, we describe a scaling for Maxwell's equations and apply it to various potential formulations and the Darwin model. Similar to \cite{raviart1996hierarchy}, we let
\beno
\bar{L} &=& \textrm{characteristic length}\\
\bar{t} &=& \textrm{characteristic time}\\
\bar{v} &=&  \textrm{characteristic speed} = \bar{L}/\bar{t}\\
\bar{\rho}, \, \bar{J} &=& \textrm{charge and current densities scaling factors}\\
\bar{E}, \,\bar{B} &=& \textrm{electric and magnetic fields scaling factors}\\
\eeno
and assume that
\begin{equation}
\label{eq:scaling}
 \epsilon_0 \frac{\bar{E}}{\bar{L}\bar{\rho}}=1, \quad \frac{\bar{J}}{\bar{\rho}}=c, \quad\frac{\bar{E}}{\bar{B}}=c.
\end{equation}
The dimensionless Maxwell's equations become
\begin{subequations}
\label{maxwells}
\begin{align}
&-\epsilon \frac{\partial  \bE}{\partial t} +\nabla \times \bB =   \bJ,\\
&\epsilon \frac{\partial  \bB}{\partial t} +\nabla \times \bE =  0,\\
&\nabla \cdot \bE =\rho,  \\
& \nabla \cdot \bB = 0,
\end{align}
\end{subequations}
where $\epsilon=\frac{\bar{v}}{c}$. The rescaled
 continuity equation is
\be
\label{continuitys}
\epsilon \frac{\partial \rho}{\partial t}+\nabla \cdot \bJ=0.
\ee
The boundary conditions are the same as \eqref{eq:pec_bc} and \eqref{eq:sl_bc} if we let $\bar{g}=\bar{E}.$
Note that the scaling considered in this paper is different from the Poiswell model \cite{masmoudi2001selfconsistent} developed for quantum mechanics, in which $-\epsilon \frac{\partial  \bE}{\partial t} +\nabla \times \bB =   \epsilon \bJ$.  Such limit is different from the regime  we are considering, and we do not consider it in this paper.
 %We also remark that the Poiswell model   developed in  for quantum mechanics is based on the scalar and vector potentials with the Lorentz gauge.

Under the same scaling \eqref{eq:scaling}, the curl-curl formulation for
the vector potential $\bw$ with
$$\bar{w}=\bar{t}\bar{E}$$
becomes
\begin{subequations}
\label{news}
\begin{align}
&\epsilon^2 \frac{\partial^2  \bw}{\partial t^2} +\nabla \times (\nabla \times \bw) =  \epsilon \nabla \times \bB(0, \bx)-\epsilon  \bJ,\\
&\nabla \cdot \left (\frac{\partial  \bw}{\partial t} \right)=\rho,\\
& \bE(t, \bx)=\partial_t \bw(t, \bx),\\
&\bB(t, \bx)=- \frac{1}{\epsilon}\nabla \times \bw(t, \bx)+\bB(0, \bx)
\end{align}
\end{subequations}
or equivalently, the wave equation formulation for $\bw$ becomes
\begin{subequations}
\label{news2}
\begin{align}
&\epsilon^2 \frac{\partial^2  \bw}{\partial t^2} -\Delta \bw = -\nabla(\int_0^t \rho ds) + \epsilon \nabla \times \bB(0, \bx)-\epsilon  \bJ,\\
&\nabla \cdot \left(\frac{\partial  \bw}{\partial t} \right)=\rho,\label{new2_continuity}\\
& \bE(t, \bx)=\partial_t \bw(t, \bx),\\
&\bB(t, \bx)=- \frac{1}{\epsilon}\nabla \times \bw(t, \bx)+\bB(0, \bx).
\end{align}
\end{subequations}
The perfectly conducting boundary condition \eqref{eq:pec_bc_w} scales as
\begin{equation}
\label{eq: bc_per_con}
\bw \times \bn=0 \quad \textrm{on} \, \,\Gamma_C,
\end{equation}
and the Silver-M\"{u}ller boundary conditions \eqref{eq:sl_bc_w} scale as
\begin{equation}
\label{eq: bc_sil_mul} (\partial_t \bw + \frac{1}{\epsilon} (\nabla \times \bw) \times \bn ) \times \bn=(\bB(0,\bx) \times \bn + \bg) \times \bn   \quad \textrm{on} \, \,\Gamma_A,
\end{equation}
if $\bar{g}=\bar{E}.$
Similarly, with the scaling
\beno
\bar{A}=\bar{B}\bar{L}, \quad \bar{\psi}=c \bar{A}=\bar{E}\bar{L},
\eeno
the evolution equations for the potentials $\psi$ and $\bA$ under the Lorentz gauge \eqref{lorentz} become
\begin{subequations}
\label{lorentzs}
\begin{align}
&\epsilon^2 \frac{\partial^2  \psi}{\partial t^2} -\Delta  \psi =  \rho,\\
&\epsilon^2 \frac{\partial^2  \bA}{\partial t^2} -\Delta  \bA =   \bJ, \\
&\nabla \cdot \bA+\epsilon\frac{\partial \psi}{\partial t}=0,\label{eq: lorentz}\\
&\bE=-\nabla \psi-\epsilon \frac{\partial \bA}{\partial t}, \quad \bB=\nabla \times \bA.
\end{align}
\end{subequations}
The evolution equations for $\psi$ and $\bA$ under the Coulomb gauge \eqref{coulombs} become
\begin{subequations}
\label{coulombs}
\begin{align}
&  -\Delta  \psi = \rho,\\
&\epsilon^2 \frac{\partial^2  \bA}{\partial t^2} -\Delta  \bA =   \bJ-\epsilon \nabla (\frac{\partial \psi}{\partial t}),\\
&\nabla \cdot \bA=0,\\
&\bE=-\nabla \psi-\epsilon \frac{\partial \bA}{\partial t}, \quad \bB=\nabla \times \bA.
\end{align}
\end{subequations}

Lastly,  the Darwin model
scales as
\begin{subequations}
\label{darwin}
\begin{align}
&-\Delta \psi =\rho, \quad \bE_{irr}=-\nabla \psi, \\
& \nabla \times \nabla \times \bB=\nabla \times \bJ, \quad \nabla \cdot \bB=0, \\
& \nabla \times \nabla \times \bE_{sol}=-\epsilon \frac{\partial}{\partial t} \nabla \times \bB, \quad \nabla \cdot \bE_{sol}=0,
\end{align}
\end{subequations}
while for the boundary conditions, we have
\beno
&&\psi|_{\Gamma_{C_i}} =\alpha_i, \quad \textrm{with} \, \,\sum_{j=0}^m c_{ij} \frac{d \alpha_j}{dt}=\frac{1}{\epsilon} \int_\Omega \bJ \cdot \nabla \chi_i \, dx, \quad 0\leq i\leq m\quad \textrm{on} \, \, \Gamma_C,  \\
&& \bB \cdot \bn= \bB(0,\bx) \cdot \bn, \quad (\nabla \times \bB) \times \bn= \bJ \times \bn \quad \textrm{on} \, \, \Gamma_C, \\
&& \bE_{sol} \times \bn=0\quad\textrm{on} \, \, \Gamma_C, \quad \int_{\Gamma_{C_i}}\bE_{sol} \cdot \bn ds=0  \quad 0\leq i\leq m,\\
&&(\nabla \times \bB- \epsilon \partial_t\bB \times \bn)\times \bn=  \bJ \times \bn +\epsilon \partial_t \bg \times \bn \quad \textrm{on} \, \, \Gamma_A,\\
&&(\epsilon \partial_t\bE_{sol}+ (\nabla \times \bE_{sol}) \times \bn)\times \bn= \epsilon \partial_t \bg \times \bn \quad \textrm{on} \, \, \Gamma_A.
\eeno

\subsubsection{Asymptotic Expansion }
\label{sec:asymptotic}
In this subsection, asymptotic expansions of the models will be performed under the scaling introduced in the previous subsection. Such expansions follow the procedure proposed in \cite{degond1992analysis}. 
For simplicity, we apply the expansions to equations in free space without boundary conditions. We impose the ansatz that variables can be expanded in terms of the dimensionless parameter $\epsilon$,
$$f=f_0+\epsilon f_1+ \epsilon^2 f_2+\cdots.$$

We first  perform an asymptotic expansion on the Darwin model \eqref{darwin}. Matching the asymptotic expansion, we obtain
\[
O(1): \quad \left\{
  \begin{array}{l}
-\Delta \psi_0=\rho_0, \\
(\bE_{irr})_0=-\nabla \psi_0, \\
 \nabla \times \nabla \times \bB_0=\nabla \times \bJ_0, \\
\nabla \cdot \bB_0=0,\\
\nabla \times \nabla \times (\bE_{sol})_0=0,\\
\nabla \cdot (\bE_{sol})_0=0,
  \end{array} \right.
  \qquad
  O(\epsilon): \quad \left\{
  \begin{array}{l}
-\Delta \psi_1=\rho_1, \\
(\bE_{irr})_1=-\nabla \psi_1, \\
 \nabla \times \nabla \times \bB_1=\nabla \times \bJ_1, \\
\nabla \cdot \bB_1=0,\\
\nabla \times \nabla \times (\bE_{sol})_1=-  \frac{\partial}{\partial t} \nabla \times \bB_0,\\
\nabla \cdot (\bE_{sol})_1=0,
  \end{array} \right.
  \]
  \[
    O(\epsilon^k), \,k \geq 2: \quad \left\{
  \begin{array}{l}
-\Delta \psi_k=\rho_k, \\
(\bE_{irr})_k=-\nabla \psi_k, \\
 \nabla \times \nabla \times \bB_k=\nabla \times \bJ_k, \\
\nabla \cdot \bB_k=0,\\
\nabla \times \nabla \times (\bE_{sol})_k=-  \frac{\partial}{\partial t} \nabla \times \bB_{k-1},\\
\nabla \cdot (\bE_{sol})_k=0.
  \end{array} \right.  \]

Similarly, we perform an asymptotic expansion for the potential formulation $\bw$ in \eqref{news} and obtain
\[
O(1): \quad \left\{
  \begin{array}{l}
\nabla \times (\nabla \times \bw_0) =0, \\
\nabla \cdot \left (\frac{\partial  \bw_0}{\partial t} \right)=\rho_0, \\
 \bE_0=\partial_t \bw_0, \\
\nabla \times \bw_0=0,\\
\nabla \cdot \bJ_0=0,
  \end{array} \right.
  \qquad
  O(\epsilon): \quad \left\{
  \begin{array}{l}
\nabla \times (\nabla \times \bw_1) =\nabla \times \bB(0, \bx)- \bJ_0, \\
\nabla \cdot \left (\frac{\partial  \bw_1}{\partial t} \right)=\rho_1, \\
 \bE_1=\partial_t \bw_1, \\
 \bB_0=- \nabla \times \bw_1+\bB(0, \bx),\\
\frac{\partial \rho_0}{\partial t}+\nabla \cdot \bJ_1=0,
  \end{array} \right.
  \]
  \[
    O(\epsilon^k), \,k \geq 2: \quad \left\{
  \begin{array}{l}
\frac{\partial^2 \bw_{k-2}}{\partial t^2}+\nabla \times (\nabla \times \bw_{k-1}) =\ - \bJ_{k-1}, \\
\nabla \cdot \left (\frac{\partial  \bw_k}{\partial t} \right)=\rho_k, \\
 \bE_k=\partial_t \bw_k, \\
 \bB_{k-1}=- \nabla \times \bw_k,\\
\frac{\partial \rho_{k-1}}{\partial t}+\nabla \cdot \bJ_k=0.
  \end{array} \right.
  \]
We now verify that the expansion for the potential formulation agrees with the Darwin model up to second order. For the $O(1)$ terms, $\nabla \cdot \bE_0=\rho_0$, and $\nabla \times \bE_0=0$. Hence, there must exist a potential function $\psi_0$, such that $\bE_0=-\nabla \psi_0$. Since $\nabla \cdot \bE_0=\rho_0$, we recover $-\Delta \psi_0=\rho_0$.  For the $O(\epsilon)$ terms, $\nabla \times \bB_0=-\nabla \times \nabla \times \bw_1+\nabla \times \bB(0, \bx)=\bJ_0$, and $\nabla \times \nabla \times \bE_1=\partial_t (\nabla \times \bB(0, \bx)-\bJ_0)=-\partial_t \bJ_0=-\partial_t \nabla \times \bB_0$, and $\nabla \cdot \bE_1=\rho_1$. Performing a Helmholtz decomposition on $\bE_1=(\bE_{irr})_1+(\bE_{sol})_1$ recovers the $O(\epsilon)$ expansion for the Darwin model. Similar derivation goes through for the wave formulation \eqref{news2} and is omitted here. %The boundary conditions could be derived similarly, and they agree with each other up to second order.

Consider now the scalar and vector potential formulations. An asymptotic expansion with the Lorentz gauge \eqref{lorentzs} gives
\[
O(1): \quad \left\{
  \begin{array}{l}
-\Delta \psi_0=\rho_0, \\
-\Delta \bA_0=\bJ_0, \\
\nabla \cdot \bA_0=0, \\
\nabla \cdot \bJ_0=0,\\
\bE_0=-\nabla \psi_0,\\
\bB_0=\nabla \times \bA_0,
  \end{array} \right.
  \qquad
  O(\epsilon): \quad \left\{
  \begin{array}{l}
-\triangle \psi_1=\rho_1, \\
-\Delta \bA_1=\bJ_1, \\
\nabla \cdot \bA_1+\frac{\partial \psi_0}{\partial t}=0, \\
\frac{\partial \rho_0}{\partial t}+\nabla \cdot \bJ_1=0,\\
\bE_1=-\nabla \psi_1-\frac{\partial \bA_0}{\partial t},\\
\bB_1=\nabla \times \bA_1.
  \end{array} \right.
  \]
  \[
    O(\epsilon^k), \,k \geq 2: \quad \left\{
  \begin{array}{l}
\frac{\partial^2 \psi_{k-2}}{\partial t^2}-\Delta \psi_k=\rho_k, \\
\frac{\partial^2 \bA_{k-2}}{\partial t^2}-\Delta \bA_k=\bJ_k, \\
\nabla \cdot \bA_k+\frac{\partial \psi_{k-1}}{\partial t}=0, \\
\frac{\partial \rho_{k-1}}{\partial t}+\nabla \cdot \bJ_k=0,\\
\bE_k=-\nabla \psi_k-\frac{\partial \bA_{k-1}}{\partial t},\\
\bB_k=\nabla \times \bA_k.
  \end{array} \right.
  \]
This expansion agrees with the Darwin model up to $O(\epsilon)$ term  if we take $(\bE_{irr})_0=\bE_0$,  $(\bE_{sol})_0=0$, $(\bE_{irr})_1=-\nabla \psi_1$, and $(\bE_{sol})_1=-\frac{\partial \bA_0}{\partial t}$. Observe
  $$ \nabla \times \bB_0=\nabla \times \nabla \times \bA_0=\nabla (\nabla \cdot \bA_0) -\Delta \bA_0=\bJ_0.$$
  Hence,
  $$ \nabla \times \nabla \times \bB_0=\nabla \times \bJ_0,$$
and
$$\nabla \cdot \bB_0=\nabla \cdot (\nabla \times \bA_0)=0$$
agree with the magneto static model on first order. However, on the first order, the Darwin equations are not consistent with the continuity equation.

Applying our procedure with the Coulomb gauge
\[
O(1): \quad \left\{
  \begin{array}{l}
-\Delta \psi_0=\rho_0, \\
-\Delta \bA_0=\bJ_0, \\
\nabla \cdot \bA_0=0, \\
\nabla \cdot \bJ_0=0,\\
\bE_0=-\nabla \psi_0,\\
\bB_0=\nabla \times \bA_0,
  \end{array} \right.
  \qquad
  O(\epsilon): \quad \left\{
  \begin{array}{l}
-\Delta \psi_1=\rho_1, \\
-\Delta \bA_1=\bJ_1-\nabla (\frac{\partial \psi_0}{\partial t}), \\
\nabla \cdot \bA_1=0, \\
\frac{\partial \rho_0}{\partial t}+\nabla \cdot \bJ_1=0,\\
\bE_1=-\nabla \psi_1-\frac{\partial \bA_0}{\partial t},\\
\bB_1=\nabla \times \bA_1,
  \end{array} \right.
  \]
  \[
    O(\epsilon^k), \,k \geq 2: \quad \left\{
  \begin{array}{l}
-\Delta \psi_k=\rho_k, \\
\frac{\partial^2 \bA_{k-2}}{\partial t^2}-\Delta \bA_k=\bJ_k-\nabla (\frac{\partial \psi_{k-1}}{\partial t}), \\
\nabla \cdot \bA_k=0, \\
\frac{\partial \rho_{k-1}}{\partial t}+\nabla \cdot \bJ_k=0,\\
\bE_k=-\nabla \psi_k-\frac{\partial \bA_{k-1}}{\partial t},\\
\bB_k=\nabla \times \bA_k.
  \end{array} \right.
  \]
  We can see that for $O(1)$ terms, the model is the same as the Lorentz gauge. For the $O(\epsilon)$ terms, a simple check yields that it still agrees with the Darwin model if we take $(\bE_{irr})_0=\bE_0$,  $(\bE_{sol})_0=0$, $(\bE_{irr})_1=-\nabla \psi_1$,  $(\bE_{sol})_1=-\frac{\partial \bA_0}{\partial t}$.

\section{Semi-discrete Schemes}
\label{sec:semidiscrete}
In this section, we extend the  implicit solver for the wave equation recently developed in \cite{causley} to Maxwell's equations. %To illustrate the main idea, we consider only the potential formulation for $\bw$ \eqref{news2} in the form of wave equations.
We  focus on the semi-discrete scheme in the MOL$^T$ framework, i.e., we only discretize the time variable and leave the space variable continuous.

\subsection{Method of Lines Transpose for Wave Equations}

The key idea of the proposed scheme is to utilize MOL$^T$ which
yields a semi-discrete system that can be solved using an integral
formulation. %The  kernel appearing in the integral formulation falls off with distance, and hence the convolution can be efficiently evaluated using a fast summation method, such as the treecode with computational complexity $\mathcal{O}(N\log N)$, where $N$ is the number of grid points (particles) used to evaluate the convolution,  see Section \ref{sec:fullydiscrete} for more details.
To illustrate the MOL$^T$ approach, consider a scalar wave equation in $\mathds{R}^3$,
 $$\frac{\partial^2u}{\partial t^2}-k^2 \Delta u =f,$$
 subject to some properly imposed boundary conditions, where $u$ is the unknown function and $k$ is a positive constant representing the wave speed.

Applying a second order finite difference approximation to $u_{tt}$ and evaluating
$\Delta u$ at time level $n+1$ gives,
$$
\Delta u^{n+1}-\frac{2}{k^2 \delta t^2} u^{n+1}= \frac{1}{k^2 \delta
  t^2 } \left ( -5 u^{n} +4u^{n-1} -u^{n-2}\right )-\frac{f^{n+1}}{k^2},
$$
with a one step truncation error of $\mathcal{O}(\delta t^3)$, where $\delta t$ denotes the time step.
$u^{n+1}$ can now be represented in an integral formulation.
\begin{align}\label{MOLT}
  u^{n+1}(\bx) &= -\int_{\Omega} \left(
      \frac{5u^n-4u^{n-1}+u^{n-2}}{k^2 \delta t^2} +\frac{f^{n+1}}{k^2}\right ) G(\bx|\by)\, d\Omega_{\by} \notag\\
  &
    -\oint_{\partial \Omega}\left( \frac{\partial u^{n+1}}{\partial\bn_\by} G(\bx|\by) - u^{n+1}\frac{\partial G}{\partial\bn_\by}(\bx|\by)\right) ds_\by,     \quad \bx\in\Omega,
\end{align}
where $G(\bx|\by)$, the free space Green's function for the modified Helmholtz
operator $\mathcal{L}(\cdot) = (\Delta - \frac{2}{k^2\delta t^2})(\cdot)$, is
$$G(\bx|\by) = -\frac{1}{4\pi r}\exp{\left(-\frac{\sqrt{2}r}{k\delta t}\right)}$$ in $\mathds{R}^3$.  Here, $r=|\bx-\by|$.  Note that, different from the fast oscillatory Green's function for the Helmholtz operator, $G(\bx|\by)$ exhibits exponential decay with respect to $r$, leading to efficient computation of the convolution integrals.  Equation \eqref{MOLT} is known as a second order dissipative scheme, which has been proposed and analyzed in \cite{causley}.

Another purely dispersive scheme with second order temporal accuracy can simply be obtained by centering the
term $\Delta u$ in time via
$$\Delta u^n\sim \frac{1}{2}\Delta (u^{n+1}+u^{n-1})~.~$$
The solution of $u^{n+1}+u^{n-1}$ is then given as,
\begin{align}
\label{MOLT2}
  (u^{n+1}+u^{n-1})(\bx) =& -\int_{\Omega}  \left(
      \frac{4u^n}{k^2 \delta t^2} +\frac{f^{n-1}+f^{n+1}}{k^2}\right ) G(\bx|\by)\, d\Omega_{\by}\\
  &- \oint_{\partial \Omega}\left( \frac{\partial (u^{n+1}+u^{n-1})}{\partial\bn_\by} G(\bx|\by) - (u^{n+1}+u^{n-1})\frac{\partial G}{\partial\bn_\by}(\bx|\by)\right) ds_\by, \quad \bx\in\Omega,\notag
\end{align}
%In this case the free space kernel needs to change ever so slightly to  $G(x,y) = -\frac{1}{4\pi|x-y|}\exp{(-\sqrt{2}|x-y|/k\Delta t)}$.

The unknown function at time level $t^{n+1}$ appearing in the boundary integrals in equations \eqref{MOLT} and \eqref{MOLT2} are then solved by imposing the boundary conditions, which will be discussed in detail in Section \ref{sec:fullydiscrete}.
Similar to \cite{christlieb2006grid} on particle-based methods for the Vlasov-Poisson system, we will evaluate the volumetric and boundary integrals using  a
midpoint  approximation.  The discrete forms of \eqref{MOLT} and \eqref{MOLT2}  can be interpreted as a collection of $N$ interacting point charges.  The resulting summations can be computed
using a fast summation algorithm, such as the  treecode algorithm to be described in Section \ref{sec:fullydiscrete}. Thanks to the implicit treatment in the MOL$^T$ approach, the proposed method is able to take time steps much larger than   an explicit integrator.

\subsection{Method of Lines Transpose for Maxwell's Equations in Potential Formulation}
\label{sec:semiasymptotic}

Using similar ideas,  we can apply the MOL$^T$ approach  to Maxwell's equations formulated using the potential formulation $\bw$.
Applying the second order dissipative MOL$^T$ approximation to \eqref{news2}, the integral solution for $\bw$ is  given by
%    \be
%    \label{scheme1ps}
%    \bw^{n+1}(x_0) &= &\iint_{\Omega} \left ( - \epsilon^2  \left(
%        \frac{2\bw^n(y)-\bw^{n-1}(y)}{ \Delta t^2} \right
%      )+T^{n+1} \right ) G(x_0|y)\, d\Omega \notag \\
%    &&+ \oint_{\partial \Omega} ((\bw^{n+1} \nabla G - G \nabla
%    \bw^{n+1}) \cdot {\bf n} \, ds
%  \ee

\begin{align}
\label{MOLTW}
  \bw^{n+1}(\bx) =& -\int_{\Omega}\left( \epsilon^2\left(
      \frac{5\bw^n-4\bw^{n-1}+\bw^{n-2}}{\delta t^2}\right) + \bT^{n+1}\right ) G(\bx|\by)\, d\Omega_{\by}\\
  &-  \oint_{\partial \Omega}\left( \frac{\partial \bw^{n+1}}{\partial\bn_\by} G(\bx|\by) - \bw^{n+1}\frac{\partial G}{\partial\bn_\by}(\bx|\by)\right) ds_\by, \quad \bx\in\Omega,\notag
\end{align}
where $\bT^n=-\nabla (\int_0^{t^n} \rho ds) +\epsilon \nabla \times \bB^0-\epsilon \bJ^{n},$ $G(\bx|\by) = -\frac{1}{4\pi r}\exp{\left(-\frac{\sqrt{2} \epsilon r}{\delta t}\right)}$ with $r=|\bx-\by|.$
If we apply the second order dispersive scheme, the solution becomes
%         \be
%    \bw^{n+1}(x_0)+ \bw^{n-1}(x_0) = \iint_{\Omega} \left ( - \epsilon^2  \left(
%        \frac{4\bw^n(y)}{ \Delta t^2} \right
%      )+\bT^{n-1}(y)+ \bT^{n+1}(y)    \right ) G(x_0|y)\, d\Omega \notag\\
%    + \oint_{\partial \Omega} ((\bw^{n+1}+\bw^{n-1}) \nabla G - G \nabla
%    (\bw^{n+1}+\bw^{n-1}) \cdot {\bf n} \, ds \label{scheme2ps}
%\ee
\begin{align}
\label{MOLT2W}
  (\bw^{n+1}+\bw^{n-1})(\bx) =& -\int_{\Omega}\left(\epsilon^2  \left(
      \frac{4\bw^n}{ \delta t^2}\right) + \bT^{n-1}+\bT^{n+1}\right ) G(\bx|\by)\, d\Omega_{\by}\\
  &- \oint_{\partial \Omega}\left( \frac{\partial (\bw^{n+1}+\bw^{n-1})}{\partial\bn_\by} G(\bx|\by) - (\bw^{n+1}+\bw^{n-1})\frac{\partial G}{\partial\bn_\by}(\bx|\by)\right) ds_\by, \quad \bx\in\Omega.\notag
\end{align}
Then, by making use of the relation between potential $\bw$ and $\bE$, $\bB$, we can further   obtain the approximations of $\bE$ and $\bB$ via $$\bE^{n+1}= \frac{3/2\bw^{n+1}-2\bw^n+\bw^{n-1}/2}{\delta t},\quad \bB^{n+1}=-\frac{1}{\epsilon} \nabla \times \bw^{n+1}+\bB^0.$$ Note that is a second order accurate scheme in time for $\bE$, and the temporal accuracy for $\bB$ is the same as $\bw$.

\begin{rem}
We mention a subtle point here for enforcing the additional constraint $\nabla \cdot \left (\frac{\partial  \bw}{\partial t} \right)=\rho.$ This is the key for charge continuity and also necessary to uniquely determine the solution as shown in Section \ref{sec:fullydiscrete}. We find that the best way is to enforce this relation using the same temporal scheme for   $\bE$, i.e. we require
\begin{equation}
\label{eq: dis_continuity}
\nabla\cdot\left(\frac{3/2\bw^{n+1}-2\bw^n+1/2\bw^{n-1}}{\delta t}\right)=\rho^{n+1}.
\end{equation}
In computations, we will only enforce \eqref{eq: dis_continuity} on the boundaries instead of the whole domain. Such practice is justified if the discrete charge and current densities $\rho, \bj$ satisfy the continuity equation (which is true if they are both zero). Otherwise, a divergence cleaning procedure is needed, and we will discuss the detailed procedure in our future work.
\end{rem}

% \textcolor{red}{This is a XX order accurate scheme for $\bE, \bB$.}

%\textcolor{red}{replace $\bB(t=0)$ in this section by $\bB^0$}

Similarly, we can  obtain the integral formulations for potential $\psi$ and $\bA$ with the Lorentz gauge and Coulomb gauge. For example, the second order dissipative scheme with the Lorentz gauge writes
\begin{align}
\label{MOLTLORENTZpsi}
  \psi^{n+1}(\bx) =& -\int_{\Omega}\left( \epsilon^2\left(
      \frac{5\psi^n-4\psi^{n-1}+\psi^{n-2}}{\delta t^2}\right) + \rho^{n+1}\right ) G(\bx|\by)\, d\Omega_{\by}\\
  &-  \oint_{\partial \Omega}\left( \frac{\partial \psi^{n+1}}{\partial\bn_\by} G(\bx|\by) - \psi^{n+1}\frac{\partial G}{\partial\bn_\by}(\bx|\by)\right) ds_\by, \quad \bx\in\Omega,\notag\\
 \label{MOLTLORENTZA}
    \bA^{n+1}(\bx) =& -\int_{\Omega}\left( \epsilon^2\left(
      \frac{5\bA^n-4\bA^{n-1}+\bA^{n-2}}{\delta t^2}\right) + \bJ^{n+1}\right ) G(\bx|\by)\, d\Omega_{\by}\\
  &-  \oint_{\partial \Omega}\left( \frac{\partial \bA^{n+1}}{\partial\bn_\by} G(\bx|\by) - \bA^{n+1}\frac{\partial G}{\partial\bn_\by}(\bx|\by)\right) ds_\by, \quad \bx\in\Omega.\notag
\end{align}
The second order dissipative scheme with the Coulomb gauge writes
\begin{align}
\label{MOLTCOULOMBpsi}
\psi^{n+1}(\bx) = &-\int_\Omega \rho^{n+1}G_0(\bx|\by)d\Omega_\by - \oint_{\partial \Omega}\left( \frac{\partial \psi^{n+1}}{\partial\bn_\by} G_0(\bx|\by) -\psi^{n+1}\frac{\partial G_0}{\partial\bn_\by}(\bx|\by)\right) ds_\by, \quad \bx\in\Omega,\\[2mm]
\label{MOLTCOULOMBA}
    \bA^{n+1}(\bx) =& -\int_{\Omega}\left( \epsilon^2\left(
      \frac{5\bA^n-4\bA^{n-1}+\bA^{n-2}}{\delta t^2}\right) + \bJ^{n+1}-\epsilon\nabla\left(\frac{3/2\psi^{n+1}-2\psi^n+\psi^{n-1}/2}{\delta t}\right)\right ) G(\bx|\by)\, d\Omega_{\by}\\
  &-  \oint_{\partial \Omega}\left( \frac{\partial \bA^{n+1}}{\partial\bn_\by} G(\bx|\by) - \bA^{n+1}\frac{\partial G}{\partial\bn_\by}(\bx|\by)\right) ds_\by, \quad \bx\in\Omega,\notag
\end{align}
where $G_0(\bx|\by)=-\frac{1}{4\pi r}$ denotes the Green's function associated with the Laplace operator $\Delta$.
Once $\psi$ and $\bA$ are obtained, we can also advance $\bE$ and $\bB$ through

$$\bE^{n+1}=-\nabla \psi^{n+1}-\epsilon \frac{3/2\bA^{n+1}-2\bA^n+\bA^{n-1}/2}{\delta t},\quad \bB^{n+1}=\nabla \times \bA^{n+1}.$$
Again, we require discrete gauge conditions, such as
$$\nabla \cdot \bA^{n+1}+\epsilon\left(\frac{3/2\psi^{n+1}-2\psi^n+1/2\psi^{n-1}}{\delta t}\right)=0,$$
for the Lorentz gauge, and
$\nabla\cdot\bA^n = 0$
for the Coulomb gauge.

\subsection{Formal Asymptotic Analysis of the Semi-discrete Schemes}
\label{sec:analysis}
%In this section, we perform the asymptotic analysis of the semi discrete schemes proposed in Section \ref{sec:semiasymptotic}.
We now verify that the semi-discrete schemes in Section \ref{sec:semiasymptotic} are AP.
 In particular,
we fix the time step size $\delta t$, and let $\epsilon \rightarrow 0$.

 We first focus on schemes \eqref{MOLTW} and \eqref{MOLT2W}.
For simplicity, we  neglect the boundary terms, but  in principle, the argument  holds when we include the boundary integrals.
 Expanding the Green's function $G(\bx|\by)$ with respect to $\epsilon$, we get
% \textcolor{red}{I think $G_1$ is wrong because of $\sqrt{2}$...}
$$G(\bx|\by)=-\frac{1}{4\pi|\bx-\by|} (1-\frac{\sqrt{2}\epsilon|\bx-\by|}{\delta t}+\ldots),$$
Therefore, $G_0(\bx|\by)=-\frac{1}{4\pi|\bx-\by|} $, $G_1(\bx|\by)=\frac{\sqrt{2}}{4\pi \delta t},\cdots$. Note that, $G_0(\bx|\by)$ is the Green's function associated with  the Laplace operator. For scheme \eqref{MOLTW}, asymptotic matching gives
$$O(1): \quad \bw_0^{n+1}(\bx)= -\int_{\Omega} \nabla \left(\int_0^{t^{n+1}} \rho_0 ds\right)     G_0(\bx|\by)\, d\Omega_{\by} $$
\begin{align*}
O(\epsilon): \quad \bw_1^{n+1}(\bx)=& -\int_{\Omega} \left(\nabla \left(\int_0^{t^{n+1}} \rho_1 ds\right)-\nabla \times \bB^0+\bJ^{n+1}_0 \right) G_0(\bx|\by) \, d\Omega_{\by} \\ &- \int_{\Omega} \nabla \left(\int_0^{t^{n+1}} \rho_0 ds\right)    G_1(\bx|\by)\, d\Omega_{\by}.
\end{align*}
Hence,
$$-\Delta \bw^{n+1}_0=\nabla\left(\int_0^{t^{n+1}} \rho_0 ds\right), \qquad -\Delta \bw^{n+1}_1=\nabla \left(\int_0^{t^{n+1}} \rho_1 ds\right)-\nabla \times \bB^0+\bJ^{n+1}_0.$$ We can also derive
$$\bE_0^{n+1}=\frac{3/2\bw_0^{n+1}-2\bw_0^n+1/2\bw_0^{n-1}}{\delta t},$$
$$\bE_1^{n+1}=\frac{3/2\bw_1^{n+1}-2\bw_1^n+1/2\bw_1^{n-1}}{\delta t},$$
$$\nabla \times \bw_0^{n+1}=0,$$
$$\bB_0^{n+1}=- \nabla \times \bw_1^{n+1}+\bB^0.$$
%\textcolor{red}{Do we have from the property of Green's function, that $\nabla \cdot  \bw^{n+1}_0=\int_0^{t^{n+1}} \rho_0 ds$}
Assume condition \eqref{eq: dis_continuity} is satisfied,
then
\begin{equation}\label{eq: w_con0}\nabla\cdot\left(\frac{3/2\bw^{n+1}_0-2\bw^n_0+1/2\bw^{n-1}_0}{\delta t}\right)=\rho_0^{n+1},\end{equation}
\begin{equation}\label{eq: w_con1}\nabla\cdot\left(\frac{3/2\bw^{n+1}_1-2\bw^n_1+1/2\bw^{n-1}_1}{\delta t}\right)=\rho_1^{n+1}.\end{equation}
This  implies
 $$\nabla\cdot\bE_0^{n+1}=\rho_0^{n+1}, \quad \nabla\cdot\bE_1^{n+1}=\rho_1^{n+1}.$$ Since $\nabla \times \bw_0^{n+1}=0,$ this gives $\nabla \times \bE_0^{n+1}=0,$ and hence there must exist $\psi_0^{n+1}$ such that $-\Delta\psi_0^{n+1} = \rho^{n+1}_0,\quad \bE^{n+1}_0 = -\nabla \psi^{n+1}_0$.
For $\bB_0^{n+1}$,
\begin{align}
\nabla\times \bB_0^{n+1}=&-\nabla\times\nabla\times \bw_1^{n+1}+\nabla\times \bB^0\notag\\
=&-\nabla\nabla\cdot \bw_1^{n+1} + \Delta \bw^{n+1}_1+\nabla\times \bB^0\notag\\
=&-\nabla\nabla\cdot \bw_1^{n+1}+ \nabla\int_0^{t^{n+1}}\rho_1 ds +\bJ_0^{n+1}\label{eq: B_0}
\end{align}
Note that \eqref{eq: w_con1} is the second order backward differentiation formulae for equation
$$\partial_t\nabla\cdot \bw_1 = \rho_1,$$
whose the exact solution is $\nabla\cdot\bw_1(t) = \int_0^{t}\rho_1ds$ due to  the zero initial condition. Hence,
$$\nabla\int_0^{t^{n+1}}\rho_1 ds-\nabla\nabla\cdot \bw_1^{n+1}  = O(\delta t^2),$$
and we obtain $$\nabla\times \bB_0^{n+1} = \bJ_0^{n+1}-\nabla\nabla\cdot \bw_1^{n+1}+ \nabla\int_0^{t^{n+1}}\rho_1 ds =  \bJ_0^{n+1} + O(\delta t^2).$$
For $\bE^{n+1}_1$,
\begin{align*}
\nabla\times\nabla\times \bE^{n+1}_1 =& \nabla\times\nabla\times\left(\frac{3/2\bw_1^{n+1}-2\bw_1^n+1/2\bw_1^{n-1}}{\delta t}\right)\\
=&\nabla\nabla\cdot\left(\frac{3/2\bw_1^{n+1}-2\bw_1^n+1/2\bw_1^{n-1}}{\delta t}\right)-\Delta\left(\frac{3/2\bw_1^{n+1}-2\bw_1^n+1/2\bw_1^{n-1}}{\delta t}\right)\\
=&\nabla\rho^{n+1}_1-\frac{1}{\delta t}\nabla\left(3/2\int_0^{t_{n+1}}\rho_1ds-2\int_0^{t_{n}}\rho_1ds+1/2\int_0^{t_{n-1}}\rho_1ds\right) \\
&-\left(\frac{3/2\bJ^{n+1}_0-2\bJ^{n}_0+1/2\bJ^{n-1}_0}{\delta t}\right).
\end{align*}
Substituting \eqref{eq: B_0} into the above equation, we obtain
\begin{align*}
\nabla\times\nabla\times \bE^{n+1}_1 =& -\nabla\times\left(\frac{3/2\bB^{n+1}_0-2\bB^{n}_0+1/2\bB^{n-1}_0}{\delta t}\right)+\nabla\rho_1^{n+1} - \nabla\nabla\cdot\left(\frac{3/2\bw^{n+1}_1-2\bw^{n}_1+1/2\bw^{n-1}_1}{\delta t}\right)\\
&=-\nabla\times\left(\frac{3/2\bB^{n+1}_0-2\bB^{n}_0+1/2\bB^{n-1}_0}{\delta t}\right)\\
&= -\partial_t\nabla\times\bB^{n+1}_0 + O(\delta t^2)
\end{align*}
due to \eqref{eq: w_con1}.
Therefore, the scheme is a consistent discretization for the Darwin model \eqref{darwin} taking into account $\nabla\cdot\bE_1^{n+1}=\rho_1^{n+1}$ and the Helmholtz decomposition for $\bE_1.$ A similar derivation holds for the dispersive scheme \eqref{MOLT2W} and is omitted.
In summary, we obtain the following theorem.
\begin{thm} If equation \eqref{eq: dis_continuity} is satisfied, the implicit schemes \eqref{MOLTW}, \eqref{MOLT2W}   or their dispersive versions will reduce to   schemes for the Darwin equation
in the limit  as  $\epsilon \rightarrow 0$. Therefore, the schemes are AP.
\end{thm}

Similarly, we can repeat the procedure for  schemes \eqref{MOLTLORENTZpsi}-\eqref{MOLTLORENTZA}, which yields
$$O(1): \qquad \psi_0^{n+1}(\bx)= \int_{\Omega} -\ \rho^{n+1}_0(\by)     G_0(\bx|\by)\, d\Omega_{\by}, $$
$$\qquad\qquad \bA_0^{n+1}(\bx)= \int_{\Omega} -\ \bJ^{n+1}_0(\by)     G_0(\bx|\by)\, d\Omega_{\by}, $$

$$O(\epsilon): \qquad \psi_1^{n+1}(\bx)= \int_{\Omega}-\rho_1^{n+1}(\by)    G_0(\bx|\by)-\rho_0^{n+1}(\by)    G_1(\bx|\by)\, d\Omega_{\by}, $$
$$\qquad\qquad \bA_1^{n+1}(\bx)= \int_{\Omega}-\bJ_1^{n+1}(\by)    G_0(\bx|\by)-\bJ_0^{n+1}(\by)    G_1(\bx|\by)\, d\Omega_{\by}.$$
Hence,
$$-\Delta \psi^{n+1}_0=\rho^{n+1}_0, \qquad -\Delta \psi^{n+1}_1=\rho_1^{n+1},$$
and
$$-\Delta \bA^{n+1}_0=\bJ^{n+1}_0, \qquad -\Delta \bA^{n+1}_1=\bJ_1^{n+1}.$$
We can also obtain
$$\bE_0^{n+1}=-\nabla \psi_0^{n+1},$$
$$\bE_1^{n+1}=-\nabla \psi_1^{n+1}- \frac{3/2\bA_0^{n+1}-2\bA_0^n+1/2\bA_0^{n-1}}{\delta t},$$
$$\bB_0^{n+1}=\nabla \times \bA_0^{n+1},$$
$$\bB_1^{n+1}=\nabla \times \bA_1^{n+1}.$$
To show that the scheme is consistent with the Darwin model, observe
\begin{align*}
\nabla \times \nabla \times \bE^{n+1}_1 =& - \nabla \times \nabla \times\left(\frac{3/2\bA_0^{n+1}-2\bA_0^n+1/2\bA_0^{n-1}}{\delta t}\right)\\
 =& - \nabla \times \left(\frac{3/2\bB_0^{n+1}-2\bB_0^n+1/2\bB_0^{n-1}}{\delta t}\right)\\
 =& - \partial_t\nabla \times\bB_0^{n+1} + O(\delta t^2).
\end{align*}
Under the assumption that the following discrete version of the Lorentz gauge condition \eqref{eq: lorentz} is satisfied, i.e.,
\begin{equation}\label{eq: dlorentz}\nabla \cdot \bA^{n+1}+\epsilon\left(\frac{3/2\psi^{n+1}-2\psi^n+1/2\psi^{n-1}}{\delta t}\right)=0,\end{equation}
then
$$\nabla \cdot \bA_0^{n+1}=0 \quad  \text{and}\quad \nabla \cdot \bA_1^{n+1}+\frac{3/2\psi_0^{n+1}-2\psi_0^n+1/2\psi_0^{n-1}}{\delta t}=0.$$
For $\bB^{n+1}_0$,
\begin{align*}
\nabla \times \bB^{n+1}_0 =&  \nabla \times \nabla \times\bA_0^{n+1} = \nabla\nabla\cdot\bA_0^{n+1} -\Delta\bA^{n+1}_0 = \bJ_0^{n+1}.
\end{align*}
Similarly, for $\bB^{n+1}_1$,
\begin{align*}
\nabla \times \bB^{n+1}_1 =&  \nabla \times \nabla \times\bA_1^{n+1} = \nabla\nabla\cdot\bA_1^{n+1} -\Delta\bA^{n+1}_1 =-\nabla\left(\frac{3/2\psi^{n+1}_0-2\psi^n_0+1/2\psi^{n-1}_0}{\delta t}\right)+ \bJ_1^{n+1}.
\end{align*}
Hence, we obtain
$$\nabla\times\nabla \times \bB^{n+1}_0 = \nabla\times\bJ_0^{n+1},\quad \nabla\times\nabla \times \bB^{n+1}_1 = \nabla\times\bJ_1^{n+1}.$$
If $(\bE_{irr})^{n+1}_0=\bE^{n+1}_0$,  $(\bE_{sol}^{n+1})_0=0$, $(\bE_{irr}^{n+1})_1=-\nabla \psi^{n+1}_1$, and $(\bE^{n+1}_{sol})_1=- \frac{3/2\bA_0^{n+1}-2\bA_0^n+1/2\bA_0^{n-1}}{\delta t},$
then
$\nabla\cdot(\bE^{n+1}_{sol})_1=- \nabla\cdot\left(\frac{3/2\bA_0^{n+1}-2\bA_0^n+1/2\bA_0^{n-1}}{\delta t}\right) = 0$, because of the gauge condition \eqref{eq: dlorentz}.
Therefore, we have shown that the semi-discrete scheme is a consistent discretization of the Darwin model up to second order. The proof for the dispersive scheme is similar and is omitted.

For the dissipative scheme with the Coulomb gauge \eqref{MOLTCOULOMBpsi}-\eqref{MOLTCOULOMBA}, we first obtain
$$-\Delta \psi^{n+1}_0=\rho^{n+1}_0, \qquad -\Delta \psi^{n+1}_1=\rho_1^{n+1},$$
$$-\Delta \bA^{n+1}_0=\bJ^{n+1}_0, \qquad -\Delta \bA^{n+1}_1=\bJ_1^{n+1}-\nabla\left(\frac{3/2\psi^{n+1}_0-2\psi^n_0+\psi^{n-1}_0/2}{\delta t}\right).$$
By performing a similar asymptotic analysis,
we can verify that, if the gauge condition
$\nabla\cdot\bA^n = 0$ holds, or $\nabla\cdot\bA^n_0=0$ and $\nabla\cdot\bA^n_1=0$, then the scheme is consistent with the Darwin model up to second order. Again, we can  prove in a similar way that the dispersive scheme is also a consistent discretization for the Darwin model.
Therefore, we arrive at the following theorem.
\medskip

\begin{thm} If  the discrete gauge conditions are satisfied, the implicit dissipative schemes \eqref{MOLTLORENTZpsi}-\eqref{MOLTLORENTZA} with the Lorentz gauge and scheme \eqref{MOLTCOULOMBpsi}-\eqref{MOLTCOULOMBA} with the Coulomb gauge
or their dispersive versions will reduce to schemes for the Darwin model
in the limit   $\epsilon \rightarrow 0$. Therefore, the schemes are AP.
\end{thm}

\section{Fully Discrete Schemes}
\label{sec:fullydiscrete}
In this section, we formulate the fully discrete schemes based on semi-discrete schemes    \eqref{MOLTW} and \eqref{MOLT2W} for the vector potential $\bw$. %Some remarks will be given to the semi-discrete scheme \eqref{MOLTLORENTZpsi}-\eqref{MOLTLORENTZA} with potential $\psi$ and $\bA$.
%The boundary integral will be evaluated by a collocation method.

\subsection{The Treecode Algorithm}
The treecode algorithm is a fast summation algorithm, which is useful for approximating the volume integrals in \eqref{MOLTW} and \eqref{MOLT2W}.
First, the domain $\Omega$ is partitioned into a set of small cells $\Omega_i$, $i=1,\ldots,\,N$. The solution $\bw$ is discretized by the associated macro particles with locations at the centers of each cell, which are denoted by $\bx_i$, $i=1,\cdots,\,N$.
%
%the solution is discretized by uniformly distributed particles $S=\{x_\bi=(x_{i_1},x_{i_2},x_{i_3})|,\,{i_1}=1,\ldots,N_1;\,{i_2}=1,\ldots,N_2;\,{i_3}=1,\ldots,N_3\}$.
Assume that the source ${\bT}$ is known at the locations of particles, which is denoted by $\bT_i$ at $\bx_i$.
 Then, the particular solution is given by the convolution of the source $\bT$ with the free space Green's function
 \begin{equation*}\Phi(\bx)=\int_\Omega \bT(\by) G(\bx|\by)d\Omega.
 \end{equation*}
The integral is discretized using the mid-point rule:
\begin{equation}
\label{eq: sum1}
\Phi(\bx)= \sum_{i}\omega_i\bT_i G(\bx|\by_i),
\end{equation}
where $\omega_i$ is the quadrature weight. $\Phi$ at the location of a particle $\bx_j$ should be calculated by
\begin{equation}
\label{eq: sum}
\Phi(\bx_j)= a\bT_j + \sum_{i\neq j}\omega_i\bT_i G(\bx_j|\by_i),
\end{equation}
where $a=\int_{\Omega_j} G(\bx_j|\by)d\Omega$. Note that even though the Green's function is singular when $\by$ approaches $\bx_j$, it is still integrable. In the simulation, integral $a$ is computed using a Gaussian quadrature rule based on the spherical coordinates. Higher order accuracy can be achieved by using a higher order quadrature formula.

Equation \eqref{eq: sum} can be computed via direct summation, which leads to computational cost on the order of $O(N^2)$. In order to speed up the calculation, a treecode \cite{li2009cartesian} can be adopted.
Below, we summarize the treecode algorithm for computing the sum
\begin{equation}
\label{eq: treecode_sum}
\phi_j=\sum_{i\neq j} q_iG(\bx_j|\bx_i),
\end{equation}
where  $q_i$ denotes the charge  associated with the particle $\bx_i$. In a treecode, particles are divided into a hierarchy of clusters. Based on the tree structure, particle-particle interactions with computational complexity of $O(N^2)$ are replaced by particle-cluster interactions with complexity of $O(N\log{N})$. Once the hierarchical clusters are formed, the treecode efficiently computes the sum \eqref{eq: treecode_sum}. For a cluster $c$ with center $\bx_c$, the contribution from $c$ can be approximated using a high-dimensional Taylor expansion.
\begin{align}
\sum_{\bx_i\in c} q_iG(\bx_j|\bx_i) &= \sum_{\bx_i\in c} q_iG(\bx_j|\bx_c + (\bx_i-\bx_c))\notag \\
&\approx \sum_{\bx_i\in c} q_i\sum_{|\bk|=0}^p \frac{1}{\bk!}\partial^\bk_{\by}G(\bx_j|\bx_c)(\bx_i-\bx_c)^\bk\notag \\
&=\sum_{|\bk|=0}^p \frac{1}{\bk!}\partial^\bk_{\by}G(\bx_j|\bx_c)\sum_{\bx_i\in c} q_i(\bx_i-\bx_c)^\bk\notag \\
&\doteq \sum_{|\bk|=0}^p  a^\bk(\bx_j,\bx_c)m_c^\bk, \label{eq: treecode_taylor}
\end{align}
where $p$ is the order of approximation, $a^\bk(\bx_j,\bx_c)$ are the Taylor coefficients of the Green's function, and $m_c^\bk$ are the cluster moments associated with the cluster $c$.  Note that the moments are independent of the particle $\bx_j$ and the Taylor coefficients are independent of the number of particles inside the cluster, leading to the speedup of the treecode. Another attractive aspect of the treecode is that the Taylor coefficients can be computed via a recurrence relation, which  significantly reduces computational cost. See \cite{li2009cartesian} for a detailed discussion.

We note that the error of the approximation \eqref{eq: treecode_taylor} is $O((r_c/R)^p)$, where $r_c=\max_{\bx_i\in c}|\bx_i-\bx_c|$ denotes the radius of the cluster and $R=|\bx_j-\bx_c|$ denotes particle-cluster distance. In particular, if $r_c/R$ is small, i.e., the cluster $c$ is considered as a far-field with respect to particle $\bx_j$, then the Taylor expansion \eqref{eq: treecode_taylor} can generate a good approximation. Otherwise, if $r_c/R$ is large, then the error becomes larger accordingly, and hence the Taylor expansion is inefficient. In the treecode, $\phi_j$ is computed using the recursive divide-conquer strategy. The standard multiple acceptance criterion (MAC) $$\frac{r_c}{R}\leq\theta$$ is adopted to determine if  cluster $c$ is a far-field, where $\theta$ is a user-specified parameter. The Taylor expansion is applied only if the MAC is satisfied, otherwise the treecode will recursively consider the interactions between particle $\bx_j$ and the children of cluster $c$. If $c$ is a leaf of the tree structure, i.e., $c$ has no children, then such a cluster is identified as a near-field and the direct summation is applied. In summary, the sum calculated via the treecode reads \cite{geng2013treecode}
\begin{equation}
\label{eq: treecode_rep}
\phi_j\approx\sum_{c\in N_j}\sum_{\bx_i\in c}q_j G(\bx_j|\bx_i) + \sum_{c\in F_j}\sum_{|\bk|=0}^p a^\bk(\bx_j,\bx_c)m_c^\bk,
\end{equation}
where $N_j$ and $F_j$ are two sets of clusters considered as near-field and far-field associated with particle $\bx_j$, respectively.
In the formulation of the proposed scheme, we also need to compute the derivatives of $\phi$, which  can be similarly obtained through the above strategy, see \cite{geng2013treecode}. For example, a high order mixed partial derivative can be calculated via
\begin{align}
\partial^\bl_\bx\phi_j& = \sum_{i\neq j} q_i\partial^\bl_\bx G(\bx_j|\bx_i)\notag\\
&\approx\sum_{c\in N_j}\sum_{\bx_i\in c} q_i \partial^\bl_\bx G(\bx_j|\bx_i) + \sum_{c\in F_j}\sum_{|\bk|=0}^p(-1)^{|\bl|}\frac{(\bk+\bl)!}{\bk!} a^{\bk+\bl}(\bx_j,\bx_c)m_c^\bk.
\label{eq: treecode_rep_der}
\end{align}

\subsection{Formulations of Fully Discrete Schemes}

In this subsection, we formulate the fully discrete schemes for  potential  $\bw$ with both direct and indirect approaches. These are common approaches   in boundary integral methods, and we refer the readers to \cite{atkinson1997numerical,sauter2011boundary} for more details.
In particular, the direct method is based on the reformulation
of the solution to \eqref{MOLTW} or \eqref{MOLT2W}, while the indirect method is based on an ansatz consisting of a single layer potential.
%Both perfectly conducting and Silver-M\"{u}ller boundary condition are discussed in the proposed framework.
We will see that both methods  can
handle  problems with the perfectly conducting boundary conditions, but the indirect method
is more convenient when dealing with the Silver-M\"{u}ller boundary conditions. We choose to illustrate the main idea of the algorithm for  the cubic domain $[0,1]^3$, but the idea can also be applied to complex geometries.
%\textcolor{red}{are the direct/indirect methods new? or what other ppl have done? is there reference?}

\subsubsection{Direct Approach for a Perfectly Conducting Cube}
We start with a simple case in which the cubic domain has six perfectly conduction faces, and $\rho$ and $\bJ$ are  0. Due to   symmetry, we only consider the first component of $\bw$ denoted by $w_1$. Since the boundary is perfectly conducting, i.e., $\bw\times \bn=0$, we can formulate a decoupled boundary condition for $w_1$. On the four boundary faces $x_2=0$, $x_2=1$, $x_3=0$, and $x_3=1$, we have $w_1=0$. In other words, on those four faces, we have the homogeneous Dirichlet boundary conditions. For the other two faces,  the continuity equation \eqref{new2_continuity} translates to the divergence-free constraint on $\bw$, i.e., $\partial_{x_1}w_1+\partial_{x_2}w_2+\partial_{x_3}w_3=0$ when $\rho=0$. Moreover, since $\bw\times\bn=0$, we have $w_2=0$ and $w_3=0$. Thus, $\partial_{x_1}w_1=0$  on the boundary faces $x_1=0$ and $x_1=1$. In other words, $w_1$ satisfies the homogeneous Neumann boundary condition $\frac{\partial w_1}{\partial \bn}=0$ on faces $x_1=0$ and $x_1=1$. Now we consider the second order dissipative scheme \eqref{MOLTW} for the $w_1$ component, i.e.,
\begin{align}
\label{eq: diss_2}
w_1^{n+1}(\bx)  =& -\int_{\Omega} \epsilon^2\left(\frac{5w_1^n-4w_1^{n-1}+w_1^{n-2}}{\delta t^2} \right) G(\bx|\by)d\Omega_\by \notag \\&- \oint_{\partial \Omega} \left (\frac{\partial w_1^{n+1}}{\partial\bn_\by} G(\bx|\by) - w_1^{n+1}\frac{\partial G}{\partial\bn_\by}(\bx|\by) \right ) ds_\by,   \quad
\bx\in\Omega
\end{align}
We denote volume integral in the formulation above as $
\phi_1(\bx)$. This term can be computed by the treecode algorithm described in the previous subsection. The unknown boundary data appearing in \eqref{eq: diss_2} are solved through the boundary integral equation. In particular, we denote the boundary faces with the homogeneous Dirichlet boundary condition by $\Gamma_D$ and the boundary faces with homogeneous Neumann boundary condition by $\Gamma_N$. We further define
the unknown Neumann trace data on $\Gamma_D$ by $\gamma_1$ and the unknown Dirichlet trace data on $\Gamma_N$ by $\gamma_2$, then \eqref{eq: diss_2} can be rewritten as
\begin{align}
%\label{eq: diss_2}
w_1^{n+1}(\bx)  =& -\phi_1(\bx) - \oint_{\partial \Omega}\left( \frac{\partial w_1^{n+1}}{\partial\bn_\by} G(\bx|\by)-w_1^{n+1}\frac{\partial G}{\partial\bn_\by}(\bx|\by)  \right)ds_\by\notag\\
= & -\phi_1(\bx) - \int_{\Gamma_D} \gamma_1(\by)G(\bx|\by)ds_\by +\int_{\Gamma_N} \gamma_2(\by)\frac{\partial G}{\partial\bn_\by}(\bx|\by) ds_\by\label{eq: representation},\quad \bx\in\Omega,
\end{align}
where the homogeneous boundary conditions have been imposed. Let $\bx$ approach the boundary, we obtain the following boundary integral equations for $\gamma_1$ and $\gamma_2$:
\begin{align}
\label{eq: diri}
0= & -\phi_1(\bx) - \int_{\Gamma_D} \gamma_1(\by)G(\bx|\by)ds_\by +\int_{\Gamma_N} \gamma_2(\by)\frac{\partial G}{\partial\bn_\by}(\bx|\by) ds_\by,\quad \bx\in\Gamma_D,\\
\frac12\gamma_2(\bx) = & -\phi_1(\bx) - \int_{\Gamma_D} \gamma_1(\by)G(\bx|\by)ds_\by  +\int_{\Gamma_N} \gamma_2(\by)\frac{\partial G}{\partial\bn_\by}(\bx|\by) ds_\by,\quad \bx\in\Gamma_N,\label{eq: neum}
\end{align}
where $\gamma_2$ in \eqref{eq: neum} is divided by 2 to account for the singular nature of double layer potential. The boundary integral equations \eqref{eq: diri} and \eqref{eq: neum} are then solved by a collocation method \cite{christlieb2006grid}. $\Gamma_D$ and $\Gamma_N$ are first divided into a set of small panels $\Gamma_{j_D}$, $j_D=1,\cdots,\,M_D$ and $\Gamma_{j_N}$, $j_N=1,\cdots,\,M_N$, respectively. The centers of panel $\Gamma_{j_D}$ and panel $\Gamma_{j_N}$ are denoted by $\bx_{j_D}$ and $\bx_{j_N}$, respectively. The unknown boundary data $\gamma_1$ and $\gamma_2$ are assumed to be constant along each panel, which are denoted by $\gamma_{1,j_D}$ and $\gamma_{2,j_N}$ on panel $\Gamma_{j_D}$ and $\Gamma_{j_N}$, respectively. In order to solve unknown boundary data $\gamma_{1,j_D}$ and $\gamma_{2,j_N}$, the boundary integral equation is discretized as
\begin{align}
\label{eq: diri_dis}
\phi_1(\bx_{i_D}) = &-\sum_{j_D=1}^{M_D} \gamma_{1,j_D}\int_{\Gamma_{j_D}}G(\bx_{i_D}|\by)ds_\by +\sum_{j_N=1}^{M_N}\gamma_{2,j_N}\int_{\Gamma_{j_N}}\frac{\partial G}{\partial\bn_\by}(\bx_{i_D}|\by) ds_\by,\quad i_D = 1,\,\cdots,\,M_D,\\
\phi_1(\bx_{i_N}) = &- \sum_{j_D=1}^{M_D}
\gamma_{1,j_D}\int_{\Gamma_{j_D}}G(\bx_{i_N}|\by)ds_\by  -\frac12\gamma_{2,i_N} +\sum_{j_N=1}^{M_N}\gamma_{2,j_N}\int_{\Gamma_{j_N}}
\frac{\partial G}{\partial\bn_\by}(\bx_{i_N}|\by) ds_\by,\quad i_N = 1,\,\cdots,\,M_N.\label{eq: neum_dis}
\end{align}
All the surface integrals  in \eqref{eq: diri_dis} and \eqref{eq: neum_dis} are approximated by a Gaussian quadrature rule. In the simulations, we use a tensor-product quadrature rule with $12^{th}$ order of accuracy. At last, we obtain a linear system for $\gamma_{1,j_D}$ and $\gamma_{2,j_N}$, which can be solved by GMRES.  Numerical experiments show that the  linear system is well-conditioned: if the error tolerance is set as $10^{-14}$, GMRES will only take 10-20 iterations to converge. The procedure for the dispersive scheme \eqref{MOLT2W} is similar and omitted. %and both methods share similar computational cost.
%\begin{align}
%w_1^{n+1}(\bx) + w_1^{n-1}(\bx)  =& - \int_{\Omega} \epsilon^2\left(\frac{4w_1^n}{\delta t^2} \right) G(\bx|\by)d\Omega_\by - \oint_{\partial \Omega}\left( \frac{\partial w_1^{n+1}}{\partial\bn_y}
% + \frac{\partial w_1^{n-1}}{\partial\bn_\by}\right) G(\bx|\by)\\&-\left( w_1^{n+1}(\by) + w_1^{n-1}(\by) \right)\frac{\partial G}{\partial\bn_\by}(\bx|\by) ds_\by.\notag
%\end{align}
%By following the similar procedure for the dissipative scheme, we can formulate a boundary integral equation for the unknown boundary data at time level $t^{n+1}$. The resulting boundary equation is solved again by a collocation method. Note that the computational cost of the two schemes is at a similar scale.

After the potential $\bw^{n+1}$ is obtained, as mentioned in Section \ref{sec:semidiscrete}, we can advance the electric field by letting $\bE^{n+1} = \frac{3/2\bw^{n+1}-2\bw^{n}+\bw^{n+1}/2}{\delta t}$. The magnetic field can be obtained by $\bB^{n+1} = -\frac{1}{\epsilon}\nabla\times \bw^{n+1} + \bB^0$. Two methods can be used to apply the $\nabla\times$ operator. We can reconstruct a local polynomial interpolating $\bw^{n+1}$, then apply the $\nabla\times$ operator to the reconstructed polynomial. In this case, the divergence-free property is attained in a discrete sense. Or we can apply the $\nabla\times$ operator to the integral representation of $\bw^{n+1}$, i.e., the operator is applied to the Green's function directly. Hence $\bB^{n+1}$ is given by an integral formulation and the divergence-free property is attained in a point-wise sense. For example, consider the second order dissipative scheme, for which the numerical solution can be written as an integral representation \eqref{eq: representation}, we obtain
\begin{align}
\label{eq: integral_B}
\bB^{n+1}(\bx) =&\bB^0+\int_\Omega\epsilon\left(\frac{5\bw^n-4\bw^{n-1}+\bw^{n-2}}{\delta t^2} \right) \nabla_\bx\times G(\bx|\by)d\Omega_\by \\&+ \frac1\epsilon\oint_{\partial \Omega}\left( \frac{\partial \bw^{n+1}}{\partial\bn_\by} \nabla_\bx\times G(\bx|\by)-\bw^{n+1}\nabla_\bx\times\frac{\partial G}{\partial\bn_\by}(\bx|\by)  \right)ds_\by.\notag
\end{align}
 Note that even though the above formulation gives an explicit representation of $\bB$, it is still subject to some numerical errors when evaluating the convolution integrals, which may result in some divergence errors at the discrete level.
In the simulations, we adopt the first method to solve for $\bB$.

So far, we have assumed that $\rho=0$. The case of  $\rho\neq0$ can be treated similarly. Note that $\bw$ should satisfy the continuity equation. Recall that, in the MOL$^T$ framework, such a constraint becomes
\begin{equation*}
        \nabla\cdot\left(\frac{3/2\bw^{n+1}-2\bw^{n}+1/2\bw^{n-1}}{\delta
        t}\right)=\rho^{n+1},
\end{equation*}
for the dissipative scheme,
or equivalently
\begin{align}
\label{eq: non_homo}
        \frac32\nabla\cdot\bw^{n+1}=\delta
        t\rho^{n+1}+2\nabla\cdot\bw^n-\frac12\nabla\cdot\bw^{n-1}
\end{align}
where   $\rho^{n+1}$ is given. We again take $w_1$ as an example.  With the perfectly conducting boundary condition,
$w_1$ still satisfies the homogeneous Dirichlet boundary condition on faces $x_2=0$, $x_2=1$, $x_3=0$, and $x_3=1$. The boundary conditions on faces $x_1=0$ and $x_1=1$ can be obtained by taking advantage of the divergence constraint \eqref{eq: non_homo}. Since $w_2$ and $w_3$ are both 0 on face $x_1=0$ and $x_1=1$, constraint \eqref{eq: non_homo} becomes
\begin{align*}
        \frac32\partial_{x_1}w_1^{n+1}=\delta
        t\,\rho^{n+1}+2\partial_{x_1}w_1^n-\frac12\partial_{x_1}w_1^{n-1},
\end{align*}
or
\begin{align*}
        \frac{\partial w_1^{n+1}}{\partial\bn}=\frac23\delta
        t\,\rho^{n+1}+\frac43\frac{\partial w_1^{n}}{\partial\bn}-\frac13\frac{\partial w_1^{n-1}}{\partial\bn},
\end{align*}
which is a non-homogeneous Neumann boundary condition on  $x_1=0$ and $x_1=1$. The  procedure can be extended similarly   to $w_2$ and $w_3$.

For complex geometries, we again assume that the boundary is discretized by a set of panels $\Gamma_j$. The unit outward normal vector $\bn_j=(n_{j,1},\,n_{j,2},\,n_{j,3})$ along panel $\Gamma_j$ is a constant. A boundary panel $\Gamma_j$ is assumed to be perfectly conducting, i.e., $\bw\times\bn_j=0$ on $\Gamma_j$, or
\begin{equation}
\label{eq: bc_per_compo}
w_2n_{j,3} - w_3n_{j,2}=0,\quad w_3n_{j,1}-w_1n_{j,3}=0,\quad w_1n_{j,2}-w_2n_{j,1}=0.
\end{equation}
For simplicity, we assume the density $\rho=0$. Hence the continuity equation will impose a divergence-free constraint on $\bw$, i.e., $\partial_{x_1}w_1+\partial_{x_2}w_2+\partial_{x_3}w_3=0.$
If $n_{j,1}=0$, $n_{j,2}$ and $n_{j,3}$ cannot be both 0, hence, $w_1=0$ on panel $\Gamma_j$ from \eqref{eq: bc_per_compo}, which is a homogeneous Dirichlet boundary condition.
If $n_{j,1}\neq0$, \eqref{eq: bc_per_compo} gives $w_2 = \frac{n_{j,2}}{n_{j,1}}w_1$ and $w_3 = \frac{n_{j,3}}{n_{j,1}}w_1$. We further take partial derivatives, and obtain $\partial_{x_2}w_2 = \frac{n_{j,2}}{n_{j,1}}\partial_{x_2}w_1$ and $\partial_{x_3} w_3 = \frac{n_{j,3}}{n_{j,1}}\partial_{x_3}w_1$, where we have used the fact that $\bn_j$ is constant along each panel. The divergence-free constraint gives
\begin{equation*}
\partial_{x_1}w_1 + \frac{n_{j,2}}{n_{j,1}}\partial_{x_2}w_1 + \frac{n_{j,3}}{n_{j,1}}\partial_{x_3}w_1 =0,
\end{equation*}
or equivalently,
$$\frac{\partial w_1}{\partial\bn_j} =0,$$
which is a homogeneous Neumann boundary condition for $w_1$. The procedure can also be applied to $w_2$ and $w_3$. In summary, on a perfectly conducting panel $\Gamma_j$, if $n_{j,k}=0$, $k=1,\,2,\, 3$, then $w_k$ satisfies the homogeneous Dirchlet boundary condition, otherwise, the homogeneous Neumann boundary condition holds true.  If $\rho\neq0$, we will get an inhomogeneous Neumann boundary condition for the case of $n_{j,k}\neq0$ instead.

\subsubsection{Indirect Approaches for a Perfectly Conducting Cube}

Now we formulate two indirect approaches for a perfectly conducting cube.
As an alternative to the representation formula \eqref{MOLTW}, we can assume  $\bw$ satisfies the following integral formulation
\begin{align}
\label{eq: indir}
\bw^{n+1}(\bx) =&-\int_{\Omega}\left( \epsilon^2\left(
      \frac{5\bw^n-4\bw^{n-1}+\bw^{n-2}}{\delta t^2}\right) + \bT^{n+1}\right ) G(\bx|\by)\, d\Omega_{\by} + \oint_{\partial\Omega} \bgamma(\by) G(\bx|\by)ds_\by\\
\doteq& -\Phi(\bx) + \oint_{\partial\Omega} \bgamma(\by) G(\bx|\by)ds_\by,\quad \bx\in\Omega,
\end{align}
where $\Phi(\bx)$ denotes the particular solution of the modified Helmholtz
equation and $\bgamma$ is an unknown vector density function associated with the
single layer potential. We now
propose two approaches to solve for $\bgamma$. The first one
is based on the fact that each component of $\bw$ should satisfy either
a Dirichlet or a Neumann boundary condition on each of the discrete boundary panels. Note
that we utilized this fact when formulating the direct method above.  We take
the first component $w_1$ as an example. %We still denote $\Gamma_{j_D}$,
%$j_D=1,\cdots,\,M_D$ and $\Gamma_{j_N}$, $j_N=1,\cdots,\,M_N$ by the sets of
%small panels where the homogeneous Dirichlet and Neumann boundary conditions are
%imposed, respectively.  The centers     are denoted by $\bx_{j_D}$ and
%$\bx_{j_N}$ accordingly. A
%similar to the direct method in the sense that Consider the perfectly conducting boundary condition, i.e., $\bw\times\bn=0$.
The discretization of the boundary and the corresponding notations are the same
as the direct method. Similar to \eqref{eq: diri}-\eqref{eq: neum}, the unknown function $\gamma$ satisfies the
following integral equations:
\begin{align}
\label{eq: in_diri}
\phi_1(\bx) = & \int_{\Gamma} \gamma_1(\by)G(\bx|\by)ds_\by , \quad \bx\in\Gamma_D,\\
\frac{\partial\phi_1(\bx)}{\partial\bn_\bx}= &  \frac12\gamma_1(\bx)+ \int_{\Gamma}
\gamma_1(\by)\frac{\partial G}{\partial\bn_\bx}(\bx|\by)ds_\by
,\quad
\bx\in\Gamma_N,\label{eq: in_neum}
\end{align}
where $\gamma_1$ is divided by 2 is to account for the singularity of the normal
derivative of the single layer potential. The corresponding discretization of
the integral equations are given by
\begin{align}
\label{eq: in_diri_dis}
\phi_1(\bx_{i_D}) = &\sum_{j=1}^{M_D}
\gamma_{1,j_D}\int_{\Gamma_{j_D}}G(\bx_{i_D}|\by)ds_\by + \sum_{j_N=1}^{M_N}
\gamma_{1,j_N}\int_{\Gamma_{j_N}}G(\bx_{i_D}|\by)ds_\by ,\quad i_D = 1,\,\cdots,\,M_D,\\
\frac{\phi_1(\bx_{i_N})}{\partial\bn_\bx} = &\frac12\gamma_{1,i_N}
+\sum_{j_D=1}^{M_D}\gamma_{1,j_D}\int_{\Gamma_{j_D}}\frac{\partial
G}{\partial\bn_\bx}(\bx_{i_N}|\by)
ds_\by+\sum_{j_N=1}^{M_N}\gamma_{1,j_N}\int_{\Gamma_{j_N}}\frac{\partial
G}{\partial\bn_\bx}(\bx_{i_N}|\by) ds_\by,\quad i_N = 1,\,\cdots,\,M_N,\label{eq: in_neum_dis}
\end{align}
where, $\phi_1(\bx_{i_D})$ and $\frac{\phi_1(\bx_{i_N})}{\partial\bn_\bx}$ are
computed via the treecode. Similar to the direct
method, the obtained linear system is well-conditioned. Once  $\bw$ is
solved, the electric field $\bE$ and magnetic field $\bB$ can be obtained through
the same procedure as described for the direct method.

The second approach  is designed to solve for the three
components of $\bw$ at the same time. Starting with the assumption
\eqref{eq: indir},
let $\bx$ approach boundary and take the cross-product with $\bn_\bx$ to obtain
\begin{align}
 \label{eq: indir_bc_per}
 0=\bw^{n+1}(\bx)\times \bn_\bx  =
-\Phi(\bx) \times\bn_\bx + \oint_{\partial\Omega} \bgamma(\by)
G(\bx|\by)\times\bn_\bx ds_\by,\quad \bx\in\Gamma_C,
 \end{align}
where the perfectly conducting boundary condition has been used. Note that $\bgamma$ is not uniquely
determined by \eqref{eq: indir_bc_per}; we still need to enforce other
constraint
to uniquely solve for $\bgamma$.  We again resort to   the divergence-free constraint for the case $\rho=0$.  The case of $\rho\neq0$ is similar and omitted. Apply the divergence operator to \eqref{eq: indir} and let $\bx$ approach boundary, we have
\begin{align}
\label{eq: indir_bc_div}
0=\nabla_\bx\cdot\bw^{n+1}(\bx)  = -\nabla_\bx\cdot\Phi(\bx) + \oint_{\partial\Omega} \bgamma(\by)\cdot \nabla_\bx G(\bx|\by)ds_\by + \frac12 \bgamma(\bx)\cdot \bn_\bx,\quad \bx\in\partial\Omega,
\end{align}
where the integral has to be understood as a Cauchy principal value. A method
proposed in \cite{guiggiani1990general} can be used to approximate the singular integral. The extra term
$\frac12 \bgamma(\bx)\cdot \bn_\bx$ accounts for the singularity of the
derivatives of the single layer potential. Combining \eqref{eq: indir_bc_per}
and \eqref{eq: indir_bc_div} gives an  integral equation for
$\bgamma$; a collocation method can be formulated
for solving the integral equations.
The numerical evidence shows that the obtained linear system is still
well-conditioned. The drawback of this approach is the need to solve a linear system with three
times larger dimensions, however, this approach can be easily extended to   the Silver-M\"{u}ller boundary
conditions as to be discussed in the next subsection.

\subsubsection{Indirect Approach for the Silver-M\"{u}ller Boundary Conditions }

The Silver-M\"{u}ller boundary conditions can be
treated as follows in a similar fashion to the perfectly conducting boundary condition for the indirect approach. In the same spirit of the MOL$^{T}$ method, we first
discretize the time variable, e.g., applying a second order finite difference
discretization to \eqref{eq: bc_sil_mul}, we obtain
\begin{equation}\label{eq: bc_sil_mul_2}
\left(\frac{3/2\bw^{n+1}-2\bw^n+1/2\bw^{n-1}}{\delta
t}+\frac{1}{\epsilon}\left(\nabla\times\bw^{n+1}\times
\bn\right)\right)\times\bn = (\bB(0,\bx)\times\bn+\bg)\times\bn,
\end{equation}
or
\begin{align}\label{eq: bc_sil_mul_3}
\left(\frac32\bw^{n+1}+\frac{\delta
t}{\epsilon}\left(\nabla\times\bw^{n+1}\times \bn\right)\right)\times\bn
&=\left(\delta t\bB(0,\bx)\times\bn+\bg+2\bw^n-\frac12\bw^{n-1}\right)\times\bn,\\
&\doteq \bR\times\bn, \notag
\end{align}
where $\bR$ is known and can be treated as a source term.  Again, we start with the ansatz \eqref{eq: indir}.
Since a single layer potential is continuous across the boundary, we have
\begin{equation}
\label{eq: bc_single}
\bw^{n+1}(\bx) = -\Phi(\bx)
+\int_{\partial\Omega}\bgamma(\by)G(\bx|\by)ds_\by,\quad\bx\in\Gamma_A.
\end{equation}
Then we apply the $\nabla_\bx\times$ operator to \eqref{eq: indir}, take the cross-product with $\bn_\bx$, and let $\bx$ approach the boundary. This gives
\begin{align}
&\nabla_\bx\times\bw^{n+1}(\bx)\times\bn_\bx = - \nabla_\bx\times\Phi(\bx)\times\bn_\bx + \frac12\bgamma(\bx) + \oint_{\partial\Omega}\nabla_\bx\times(\bgamma(\by) G(\bx|\by))\times\bn_\bx ds_\by\notag\\
\label{eq: bc_curl}
& = - \nabla_\bx\times\Phi(\bx)\times\bn_\bx + \frac12\bgamma(\bx) +
\oint_{\partial\Omega} \bgamma(\by)\frac{\partial G}{\partial\bn_\bx}(\bx|\by) -
\bgamma(\by)\cdot\bn_\bx\nabla_\bx G(\bx|\by) ds_\by,\quad \bx\in\Gamma_A,
\end{align}
where  the extra term $\frac12\bgamma(\bx)$ is due to the singularity of the derivative of the single layer potential. Again the integral exists as a Cauchy principal value.
Substituting \eqref{eq: bc_single} and \eqref{eq: bc_curl} into  \eqref{eq: bc_sil_mul_3} gives
\begin{align}
\label{eq: bc_sil_mul_int_eqn}
&\left(\frac32\Phi(\bx) + \frac{\delta t}{\epsilon}\nabla_\bx\times\Phi(\bx)\times\bn_\bx+ \bR\right)\times\bn_\bx \notag\\
=&\left(\frac{\delta t}{2\epsilon}\bgamma(\bx)+ \oint_{\partial\Omega}\left(\frac{3}{2}\bgamma(\by)G(\bx|\by)+
\frac{\delta t}{\epsilon}\left(\bgamma(\by)\frac{\partial G}{\partial\bn_\bx}(\bx|\by)-\bgamma(\by)\cdot\bn_\bx\nabla_\bx G(\bx|\by)\right)\right)\,ds_\by\right)\times\bn_\bx, \quad \bx\in\Gamma_A.
\end{align}
Similar to the case of the perfectly conducting boundary condition, \eqref{eq:
bc_sil_mul_int_eqn} cannot uniquely determine $\bgamma$;
the continuity equation on the boundary $\Gamma_A$, i.e. \eqref{eq: indir_bc_div} should be enforced for the case when $\rho=0$. We remark that it is not trivial to
formulate a proper boundary integral equation in the setting of the direct
method or the first approach of the indirect method, when the
Silver-M\"{u}ller boundary condition is considered. Also note that, unlike the perfectly conducting boundary condition, the time derivative $\partial_t\bw$ appears in the Silver-M\"{u}ller boundary condition, which makes it difficult to formulate a dispersive scheme. For simplicity, we only consider the dissipative scheme for the problem with the Silver-M\"{u}ller boundary condition.

Lastly, when both types of boundary conditions are imposed,
we have \begin{align}
 \label{eq: bc_mix_per}
&\Phi(\bx) \times\bn_\bx = \oint_{\partial\Omega} \bgamma(\by) G(\bx|\by)\times\bn_\bx ds_\by, \quad \bx \in\Gamma_C,\\
&\left(\frac32\Phi(\bx) + \frac{\delta t}{\epsilon}\nabla_\bx\times\Phi(\bx)\times\bn_\bx+ \bR\right)\times\bn_\bx \notag\\
=&\left(\frac{\delta t}{2\epsilon}\bgamma(\bx)+\oint_{\partial\Omega}\frac{3}{2}\bgamma(\by)G(\bx|\by)+
\frac{\delta t}{\epsilon}\left(\bgamma(\by)\frac{\partial G}{\partial\bn_\bx}(\bx|\by)-\bgamma(\by)\cdot\bn_\bx\nabla_\bx G(\bx|\by)\right)\,ds_\by\right)\times\bn_\bx, \quad\bx\in\Gamma_A,\label{eq: bc_mix_sil_mul}
\end{align}
together with the divergence constraint on $\bw$ \eqref{eq: indir_bc_div}. The resulting system can be computed by collocation methods and linear solvers. This completes the description of our algorithms.

\subsubsection{Remarks on the Formulations with $\psi$ and $\bA$}

At the end of this section, we remark on the fully discrete schemes in terms of potential $\psi$ and $\bA$. Since we can not obtain decoupled boundary conditions for $\psi$ and $\bA$, it is impossible to formulate a direct approach; the only practical way is to  adopt the indirect approach. Moreover, similar to $\bw$, in order to uniquely determine $\psi$ and $\bA$, we need to enforce the gauge condition in the integral formulations, which will result in a linear system with dimension one third  larger than that from the formulation of $\bw$. Therefore, the use of $\bw$ seems more efficient than $\psi$ and $\bA$. On the other hand, the time integral $\nabla(\int_0^t \rho ds)$ in the formulation of $\bw$ may lead to potential difficulty for  plasma simulations.

\section{Numerical Examples}
\label{sec:numerical}

In this section, we consider two numerical examples to demonstrate the performance of the proposed
schemes.  The rescaled Maxwell's equations \eqref{maxwells} are solved numerically on   a  unit cube $[0,1]^3$ with different initial and boundary conditions.

\textbf{Problem 1.} The first test case \cite{ciarlet2007continuous} has    perfectly conducting boundaries on all six faces. The charge and current densities are set to be zero. The initial conditions for $\bE$ and $\bB$ are given by
\begin{equation*}
\bE(0,\bx)=\left(\begin{array}{c}\cos(\pi x_1)\sin(\pi x_2)\sin(-2\pi x_3)\\
\sin(\pi x_1)\cos(\pi x_2)\sin(-2\pi x_3)\\
\sin(\pi x_1)\sin(\pi x_2)\cos(-2\pi x_3)
 \end{array}\right),\quad \bB(0,\bx) = 0.
\end{equation*}
Note that both $\bE$ and $\bB$ satisfy the divergence-free condition initially. The exact solution is
\begin{align*}
\bE(t,\bx)&=\cos(\omega t)\bE(0,\bx),\\
\bB(t,\bx)&=\sqrt{\frac{3}{2}}\sin(\omega t)\left(\begin{array}{c}-\sin(\pi x_1)\cos(\pi x_2)\cos(-2\pi x_3)\\
\cos(\pi x_1)\sin(\pi x_2)\cos(-2\pi x_3)\\
0
\end{array}\right),
\end{align*}
%\begin{align*}
%\bE(t,\bx)&=\cos(\omega t)\left(\begin{array}{c}\cos(\pi x_1)\sin(\pi x_2)\sin(-2\pi x_3)\\
%\sin(\pi x_1)\cos(\pi x_2)\sin(-2\pi x_3)\\
%\sin(\pi x_1)\sin(\pi x_2)\cos(-2\pi x_3)
% \end{array}\right),\\
%\bB(t,\bx)&=\sqrt{\frac{3}{2}}\sin(\omega t)\left(\begin{array}{c}-\sin(\pi x_1)\cos(\pi x_2)\cos(-2\pi x_3)\\
%\cos(\pi x_1)\sin(\pi x_2)\cos(-2\pi x_3)\\
%0
% \end{array}\right),
%\end{align*}
where $\omega = \frac{\sqrt{6}\pi}{\epsilon}$. Recall that $\bw(t,\bx)=\int_0^t \bE(s,\bx)dt$. The exact solution for $\bw$ can be obtained by integrating $\bE$ in time:
%\begin{align*}
%\bw(t,\bx)&=\frac{1}{\omega}\sin(\omega t)\left(\begin{array}{c}\cos(\pi x_1)\sin(\pi x_2)\sin(-2\pi x_3)\\
%\sin(\pi x_1)\cos(\pi x_2)\sin(-2\pi x_3)\\
%\sin(\pi x_1)\sin(\pi x_2)\cos(-2\pi x_3)
% \end{array}\right).
% \end{align*}
%
\begin{align*}
\bw(t,\bx)&=\frac{1}{\omega}\sin(\omega t)\bE(0,\bx).
\end{align*}

\textbf{Problem 2.} The second test case \cite{assous1996numerical} is a cubic waveguide, in which a TE$_{10}$ mode propagates in the $x_3$-direction. The charge and current densities are set to be zero.  The analytical expression of the TE field is given by
\begin{align*}
\bE(t,\bx)=\left(\begin{array}{c}\sin(\pi x_2)\sin(\pi x_3-\omega t)\\
0\\
0
 \end{array}\right),\quad
\bB(t,\bx)=\frac{1}{\sqrt{2}}\left(\begin{array}{c}0\\
\sin(\pi x_2)\sin(\pi x_3-\omega t)\\
\cos(\pi x_2)\cos(\pi x_3-\omega t)
 \end{array}\right),\\
\end{align*}
where $\omega = \frac{\sqrt{2}\pi}{\epsilon}$. The exact expression for $\bw$ can be obtained by integrating $\bE(t)$ in time
\begin{align*}
\bw(t,\bx)=\left(\begin{array}{c}\frac1\omega\sin(\pi x_2)(\cos(\pi x_3-\omega t)-\cos(\pi x_3))\\
0\\
0
 \end{array}\right).
\end{align*}
We prescribe the perfectly conducting boundary condition on the four side faces ($x_1=0$, $x_1=1$, $x_2=0$, and $x_2=1$), while at the bottom ($x_3=0$) and  the top ($x_3=1$),  the Silver-M\"{u}ller boundary conditions are imposed, i.e.,
\begin{align*}
&\left(\partial_t\bw +\frac{1}{\epsilon}\nabla\times\bw\times\bn\right)\times\bn = (\bB(0,\bx)\times\bn+ \bg)\times\bn,\quad\text{with}\\
&\bg =\left\{\begin{array}{c} \displaystyle (-\frac{\sqrt{2}+2}{2}\sin(\pi x_2)\sin(\omega t),\,0,\,\,0)^T\quad\text{on}\quad x_3=0,\\[2mm]
\displaystyle  (-\frac{\sqrt{2}-2}{2}\sin(\pi x_2)\sin(\omega t),\,0,\,\,0)^T\quad\text{on}\quad x_3=1.\end{array}\right.
\end{align*}

\textbf{Problem 3.} The last test case considered in this paper is a simple case of the problem 2 studied in \cite{ciarlet2007continuous}. The computational domain is still a unit cube, but a current bar with zero thickness is posed  across the domain in the $x_3$-direction. We let $\bJ=(0,0,\delta(x_1-1/2)\delta(x_2-1/2)\cos(2\pi t))$ and $\rho=0$. We prescribe the absorbing Silver-M\"{u}ller boundary conditions on the four side faces ($x_1=0$, $x_1=1$, $x_2=0$, and $x_2=1$), i.e., we let $g=0$, while at the bottom ($x_3=0$) and  the top ($x_3=1$),  the perfectly conducting boundary condition is imposed. We set zero initial conditions for both $\bE$ and $\bB$, and hence $\bw$ is set to be zero.

\bigskip
In  our simulations, we consider both the direct
 and indirect approaches. For simplicity, we use uniformly distributed particles in the unit cube and let $h$   be the mesh size in one coordinate direction.
However, the scheme can also accommodate non-uniformly distributed particles.
 	 The time step $\delta t$ is chosen as $\delta t = \text{CFL} \cdot h.$
Note that for the test cases, choosing a small dimensionless parameter $\epsilon$ is equivalent to choosing a large CFL number, since the frequency $\omega=O(\frac{1}{\epsilon})$. Therefore, we set $\epsilon=1$ in all simulations and test the   schemes with large CFL numbers. For the treecode algorithm, the MAC parameter $\theta$ is set to be 0.5 and the order of Taylor approximation $p$ is set to be 9.

\subsection{Direct Approach for Problem 1}
\begin{figure}[!hptb]
\begin{center}
\includegraphics[width=3in]{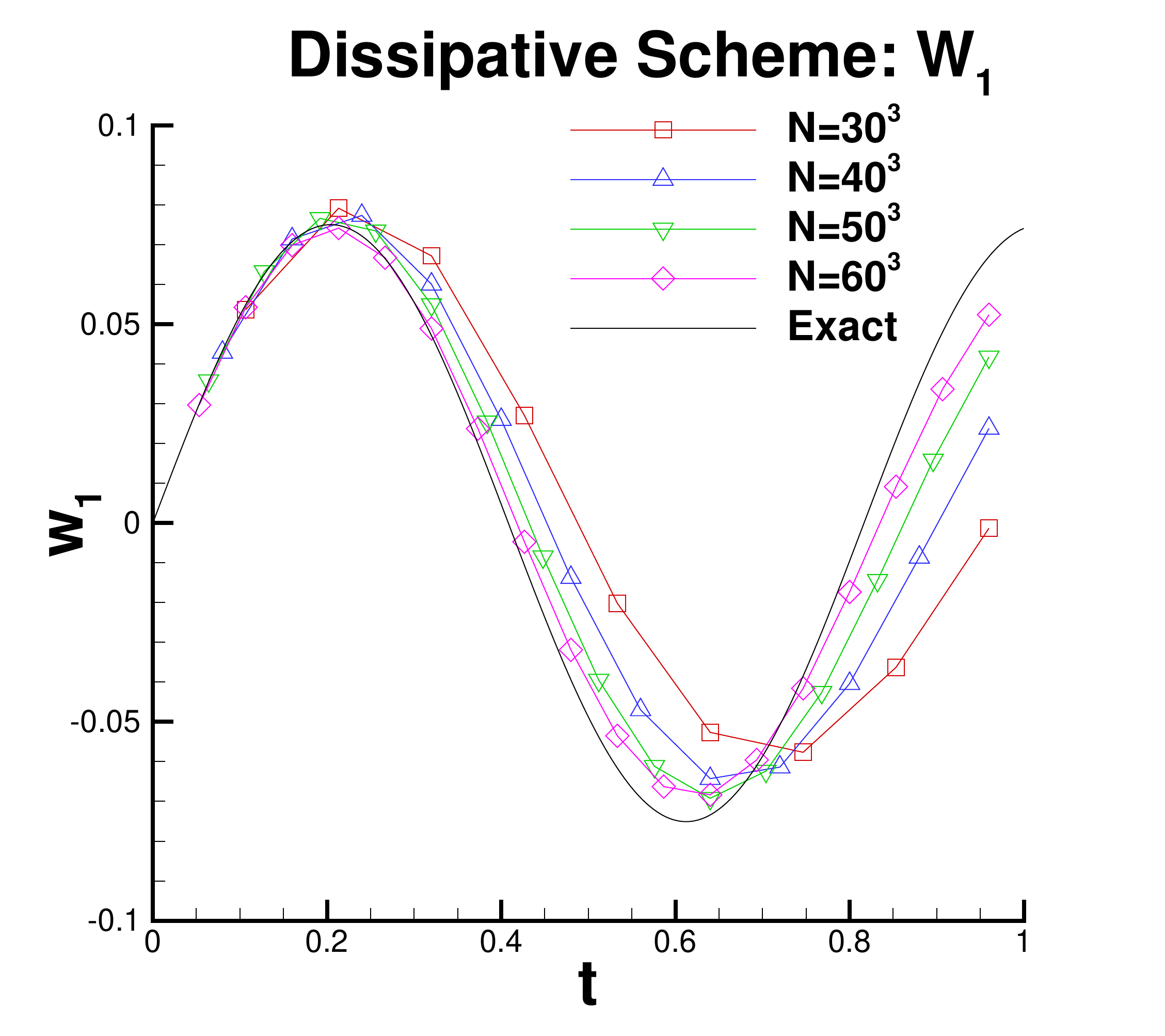}
\includegraphics[width=3in]{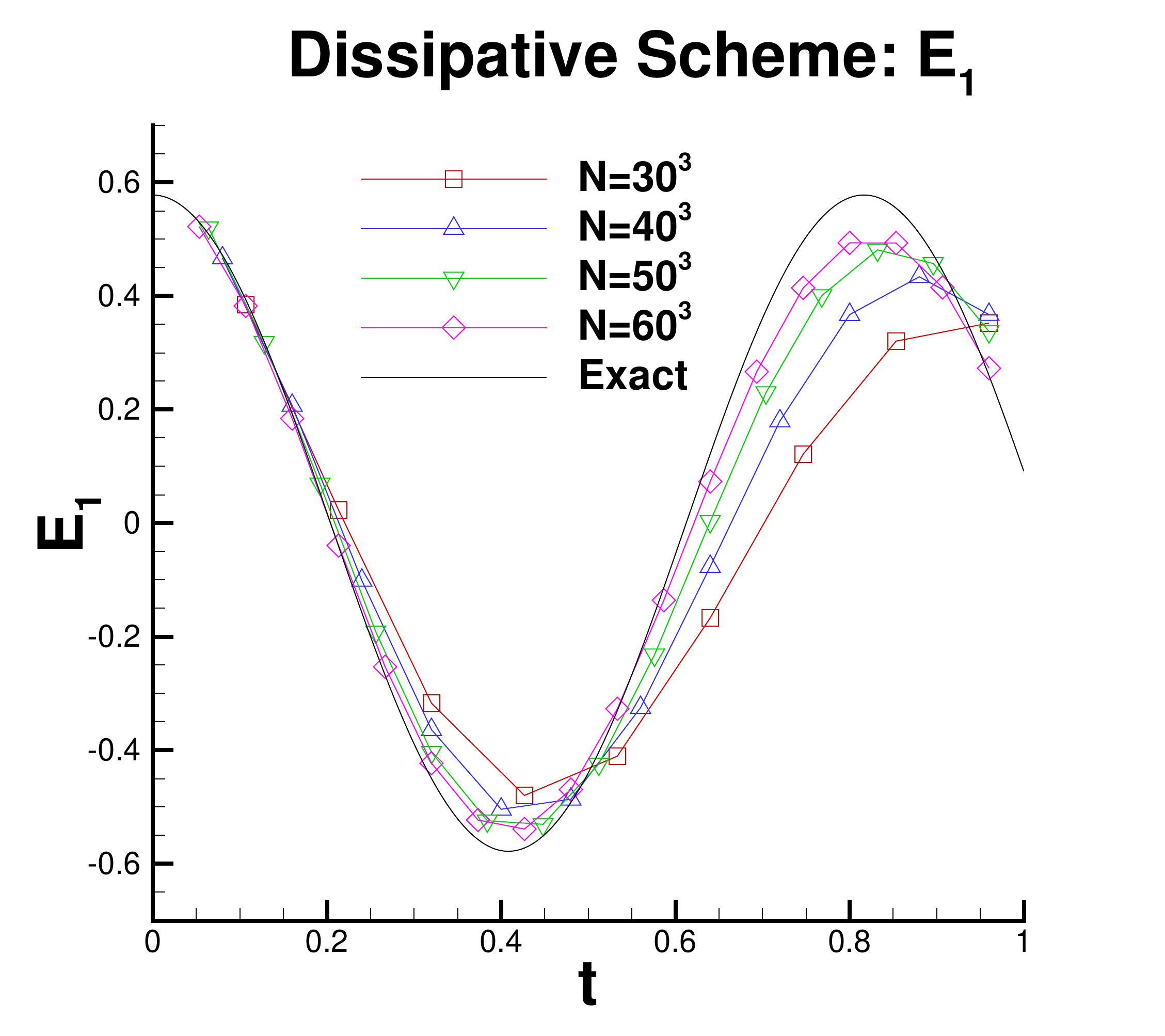}
\end{center}
\caption{
Problem 1.  The time evolution of $w_1$ (left) and $E_1$ (right) at the location ($\sqrt{3}/2,\,\sqrt{2}/2,\,\sqrt{2}/4$). Second order dissipative scheme. Direct approach. $\text{CFL} =3.2$.
}
\label{fig:time_evo_diss_dir}
\end{figure}

For the first numerical test, we let CFL=3.2 for the dissipative scheme. In Figure \ref{fig:time_evo_diss_dir}, we plot the time evolution of the numerical solutions $w_1$ and $E_1$ at an arbitrarily chosen point ($\sqrt{3}/2,\,\sqrt{2}/2,\,\sqrt{2}/4$) computed with the different number of particles. The solution value at this point is calculated from the integral formulation \eqref{eq: diri_dis}-\eqref{eq: neum_dis} once we solve the unknown boundary potential. 
It is observed the numerical solution converges to the exact solution when adding more particles. Due to the symmetry of this test, we only report the first component of the numerical solutions for brevity. In Table \ref{tab:diss_dir}, we report the convergence study of the proposed scheme. Second order of convergence for $w_1$, $E_1$ and $B_1$ are observed, where errors are measured using the following norm
 \begin{equation}
 \|e(\bx,t)\| = \|(\|e(\bx,t)\|_{L^1(\Omega)})\|_{L^\infty[0,T]}.
 \end{equation}
 %\textcolor{red}{This is L2 error, but the tables say L1 error?}

\begin{table}[!hptb]
%\label{tab: tab_burgers_1}
\begin{center}
\caption{Problem 1.   The $L^1$ errors and the corresponding orders of accuracy for $w_1$, $E_1$, and $B_1$. $T=1$.  Second order dissipative scheme. Direct approach. CFL=3.2.
\label{tab:diss_dir}}
\bigskip
\begin{tabular}{|c | c c|c c|c c|}
\hline
\cline{2-7} &\multicolumn{2}{c|}{$w_1$} &\multicolumn{2}{c|}{$E_1$} &\multicolumn{2}{c|}{$B_1$}\\
\cline{2-7}N & $L^1$ error&order&  $L^1$ error&order&  $L^1$ error&order\\
\hline
$30^3$	&3.06E-02&	--   & 	1.68E-01&	--	&2.71E-01&	--	\\ \hline
$40^3$	&1.96E-02&	1.55&	1.10E-01&	1.47&	1.78E-01&	1.46\\ \hline
$50^3$	&1.23E-02&	2.08&	6.65E-02&	2.25&	1.13E-01&	2.03\\ \hline
$60^3$	&8.33E-03&	2.14&	4.33E-02&	2.35&	7.48E-02&	2.33
\\ \hline

\hline
\end{tabular}
\end{center}
\end{table}

\begin{figure}[!hptb]
\begin{center}
\includegraphics[width=3in]{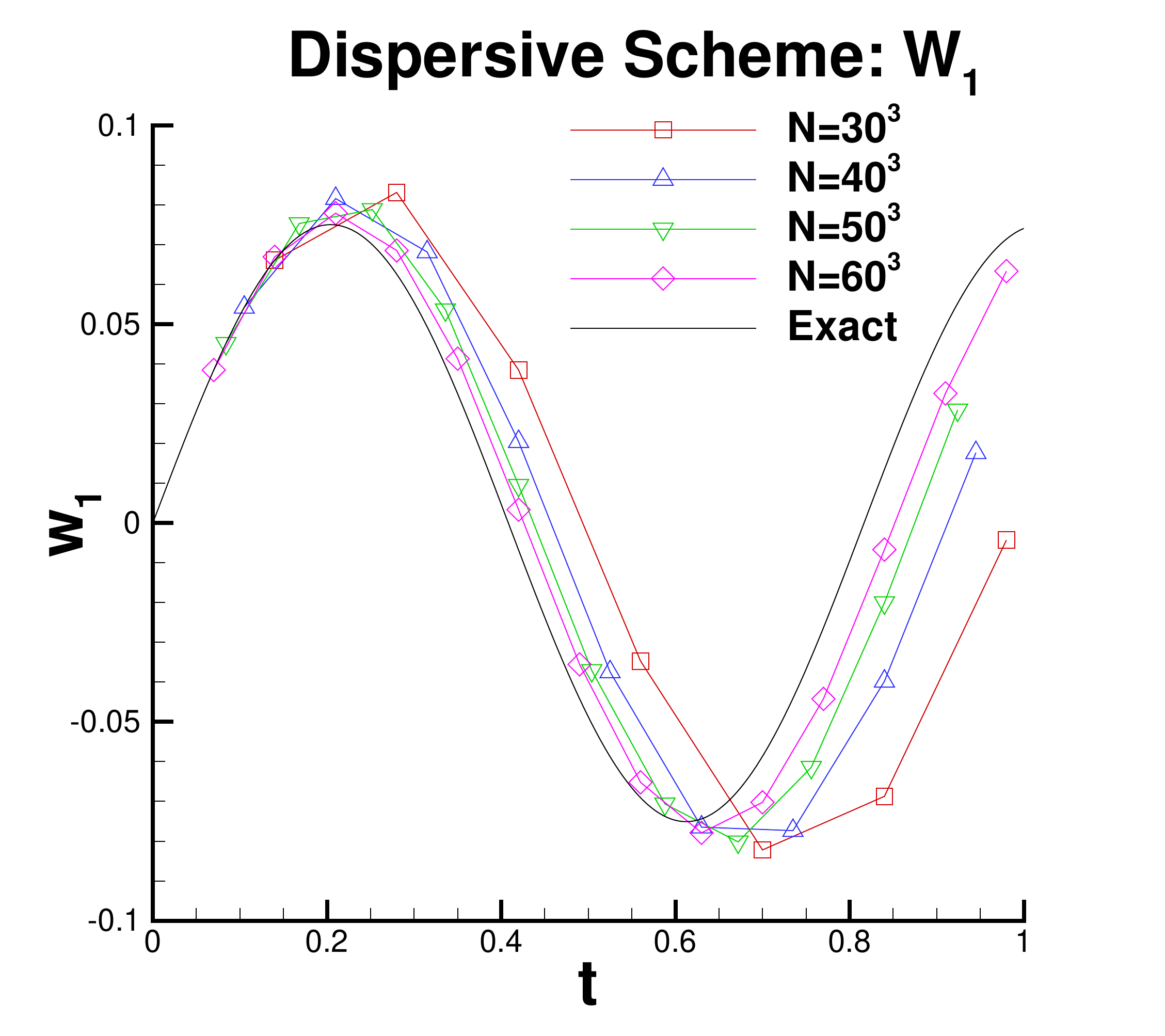}
\includegraphics[width=3in]{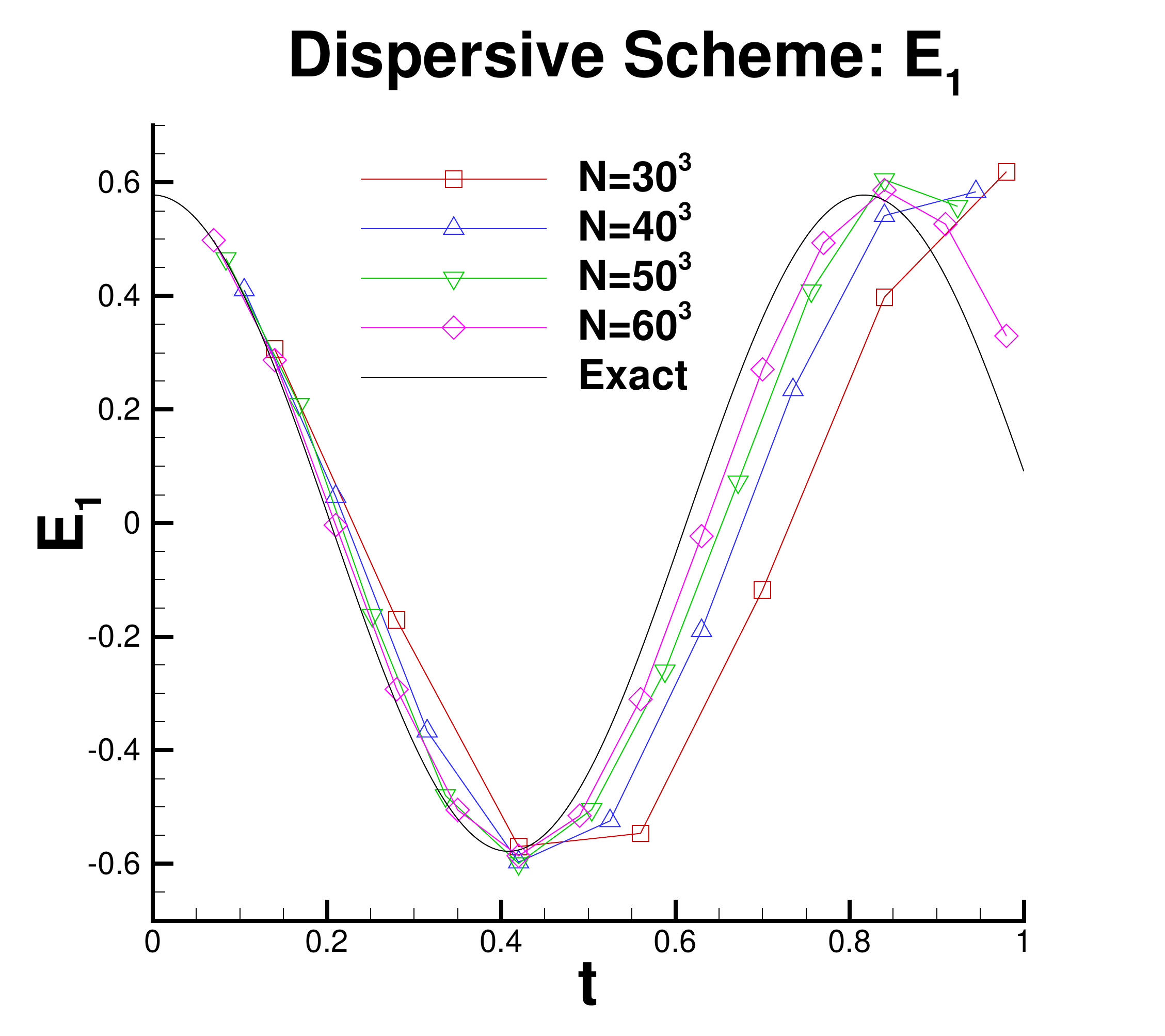}
\end{center}
\caption{
Problem 1.  The time evolution of $w_1$ (left) and $E_1$ (right) at the location ($\sqrt{3}/2,\,\sqrt{2}/2,\,\sqrt{2}/4$). Second order dispersive scheme.   Direct approach.  $\text{CFL} =4.2$.
}
\label{fig:time_evo_disp_dir}
\end{figure}

\begin{figure}[!hptb]
\begin{center}
\includegraphics[width=3in]{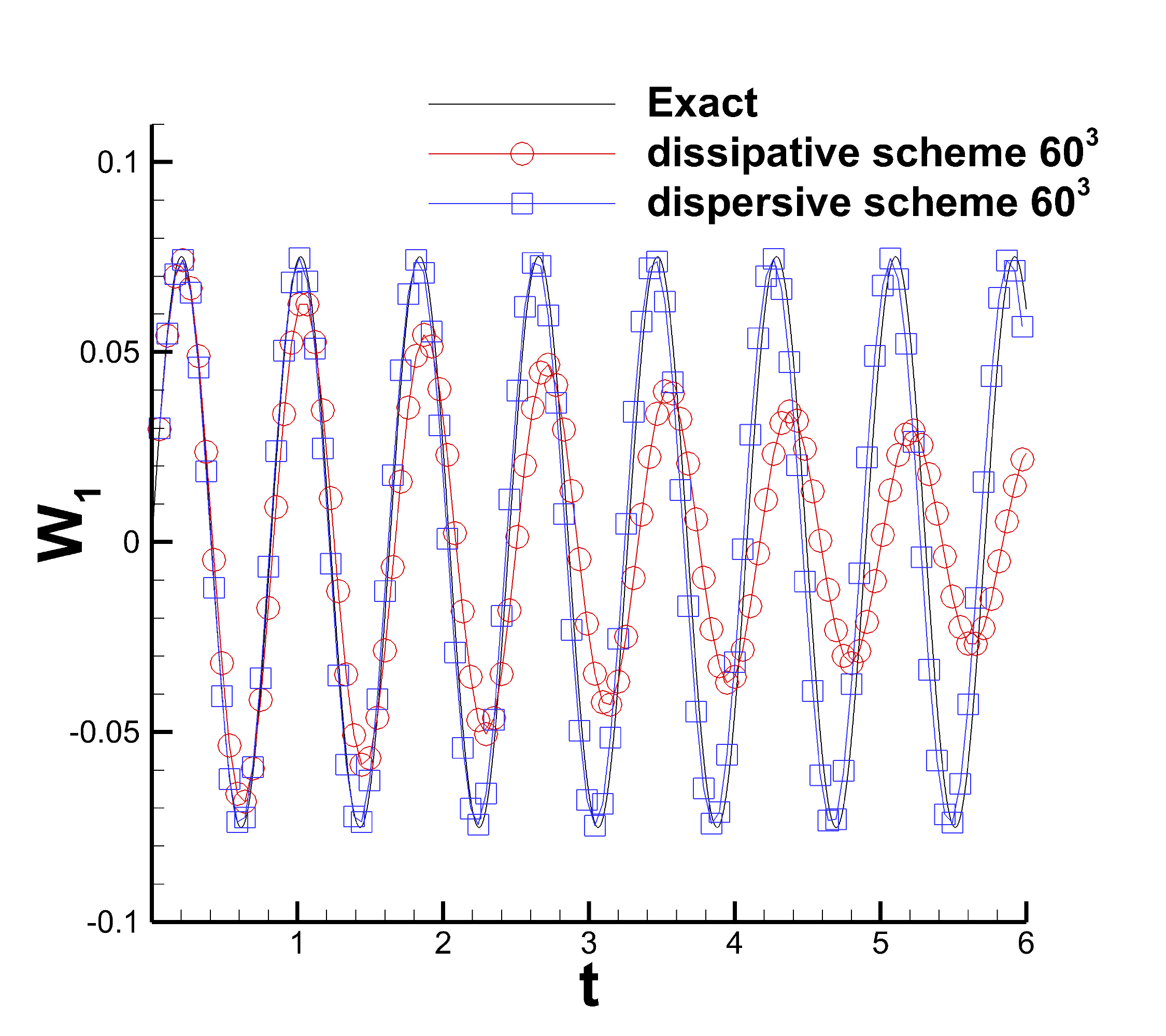}
\includegraphics[width=3in]{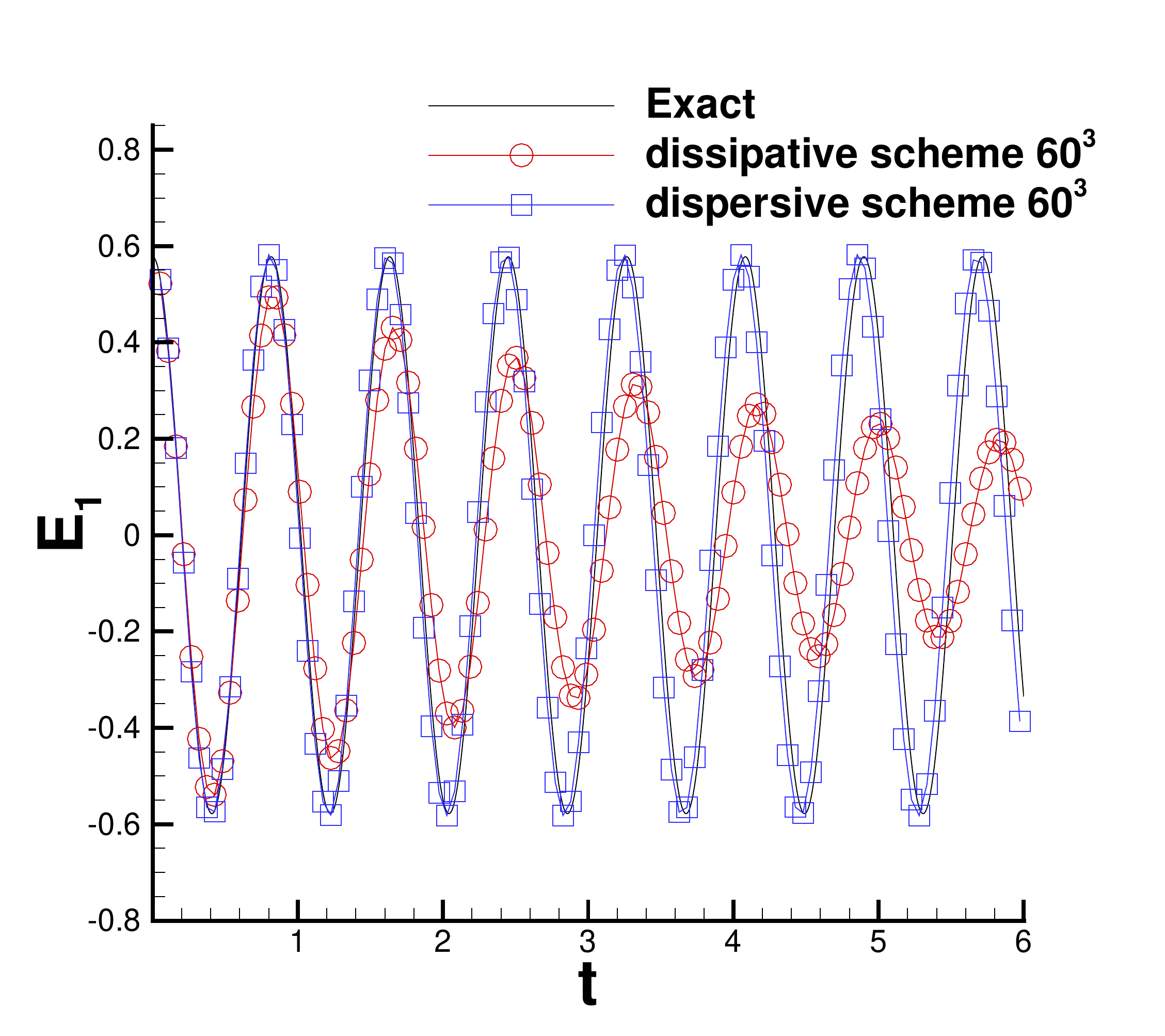}
\end{center}
\caption{
Problem 1.   The time evolution of $w_1$ (left) and $E_1$ (right) at the location ($\sqrt{3}/2,\,\sqrt{2}/2,\,\sqrt{2}/4$) up to about 7 periods. $60^3$ particles are used for computation. Direct approach. $\text{CFL} =3.2$.
}
\label{fig:time_evo}
\end{figure}

We also ran the simulations using the second order dispersive scheme for this test. Note that, we can use a larger CFL number for the dispersive scheme and obtain
comparable numerical results. In the simulation, we set CFL=4.2. In Figure \ref{fig:time_evo_disp_dir}, the time evolution of the numerical solution at the location ($\sqrt{3}/2,\,\sqrt{2}/2,\,\sqrt{2}/4$) is reported. Again, using more particles can generate
more accurate solutions. The convergence study is also presented in Table \ref{tab:disp_dir}, and second order of convergence is observed as expected. Then, we compare the performance of two
second order schemes for a long time simulation, for which we compute the numerical solution up to 7 periods. We set CFL=3.2 and use $60^3$ particles for both schemes. In Figure \ref{fig:time_evo},   we plot the numerical solutions $w_1$ and $E_1$ by both schemes at the location ($\sqrt{3}/2,\,\sqrt{2}/2,\,\sqrt{2}/4$). It is observed that the amplitude of the wave is dissipated for the dissipative scheme and the corresponding dissipation error becomes significant after a long time simulation. On the other hand, the dispersive scheme can maintain the amplitude of the wave to some extent, while the phase error can be observed after some time.

\begin{table}[htb]
%\label{tab: tab_burgers_1}
\begin{center}
\caption{Problem 1. The $L^1$ errors and the corresponding orders of accuracy for $w_1$, $E_1$, and $B_1$.  $T=1$.  Second order dispersive scheme. Direct approach.  CFL=4.2.
\label{tab:disp_dir}}
\bigskip
\begin{tabular}{|c | c c|c c|c c|}
\hline
\cline{2-7} &\multicolumn{2}{c|}{$w_1$} &\multicolumn{2}{c|}{$E_1$} &\multicolumn{2}{c|}{$B_1$}\\
\cline{2-7}N & $L^1$ error&order&  $L^1$ error&order&  $L^1$ error&order\\
\hline
$30^3$	&3.68E-02&	--   & 	2.14E-01&	--	&3.27E-01&	--	\\ \hline
$40^3$	&2.37E-02&	1.53&	1.21E-01&	1.98&	2.16E-01&	1.44\\ \hline
$50^3$	&1.48E-02&	2.11&	7.41E-02&	2.19&	1.36E-01&	2.07\\ \hline
$60^3$	&8.99E-03&	2.73&	4.52E-02&	2.71&	8.31E-02&	2.70\\
\hline
\end{tabular}
\end{center}
\end{table}

\subsection{Indirect Approach for Problem 1}

\begin{figure}[!hptb]
\begin{center}
\includegraphics[width=3in]{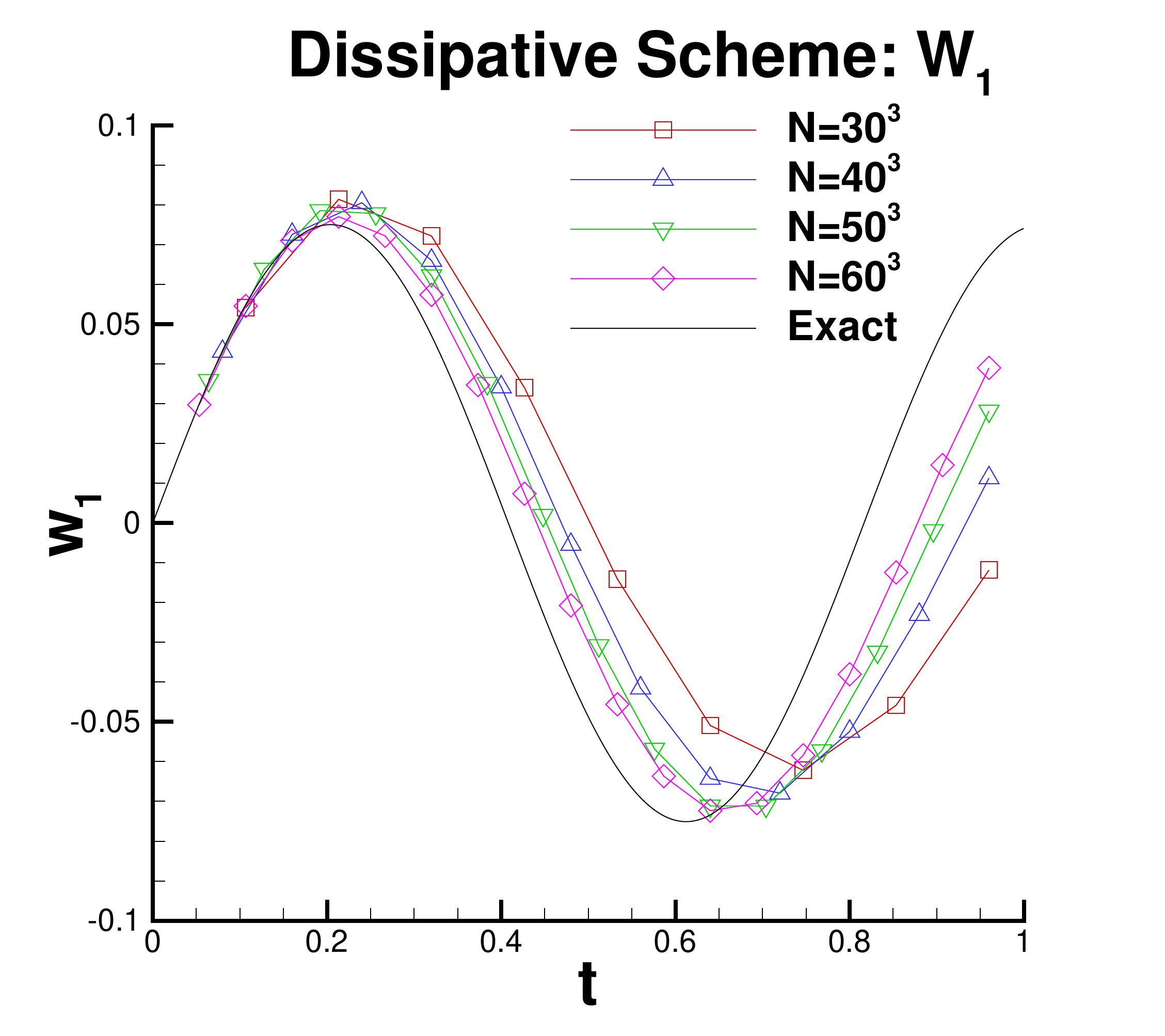}
\includegraphics[width=3in]{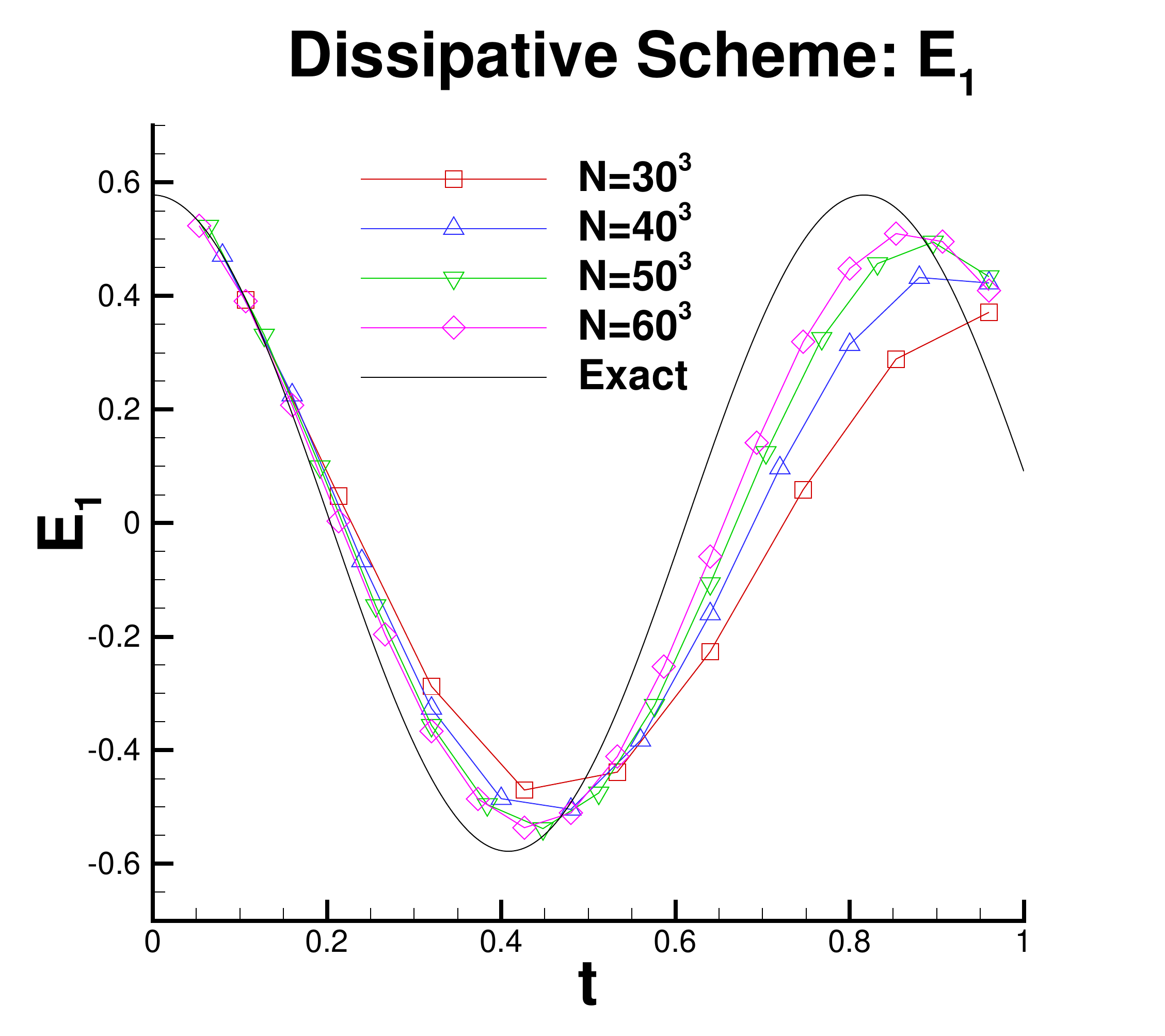}
\end{center}
\caption{
Problem 1.  The time evolution of $w_1$ (left) and $E_1$ (right) at the location ($\sqrt{3}/2,\,\sqrt{2}/2,\,\sqrt{2}/4$).  Second order dissipative scheme.  The first indirect approach. $\text{CFL} =3.2$.
}
\label{fig:time_evo_in}
\end{figure}

 We use the first indirect method to solve the first problem and set CFL=3.2.  For brevity, we only consider the second order dissipative scheme. In Figure \ref{fig:time_evo_in}, we plot the time evolution of numerical solutions $w_1$ and $E_1$ at the location  ($\sqrt{3}/2,\,\sqrt{2}/2,\,\sqrt{2}/4$). Comparable numerical results are observed to the direct method. In Table \ref{tab:diss_in}, we report the convergence study for the indirect scheme, for which we observe $1.5^{th}$ order of convergence for $w_1$ and $E_1$, and second order of convergence for $B_1$. We will investigate the reason for the reduction of accuracy in the future. We also noted that the magnitude of errors is a little larger than that by the direct method. For brevity, we do not report the numerical result for the second indirect approach, which gives comparable numerical results to the first indirect approach. %As we mentioned in Section \ref{sec:fullydiscrete}, the second approach is more costly due to its three times larger dimensions of the linear system.

\begin{table}[!hptb]
%\label{tab: tab_burgers_1}
\begin{center}
\caption{ Problem 1.    The $L^1$ errors and the corresponding orders of accuracy for $w_1$, $E_1$, and $B_1$. $T=1$.  Second order dissipative scheme.  The first indirect approach. CFL=3.2.
\label{tab:diss_in}}
\bigskip
\begin{tabular}{|c | c c|c c|c c|}
\hline
\cline{2-7} &\multicolumn{2}{c|}{$w_1$} &\multicolumn{2}{c|}{$E_1$} &\multicolumn{2}{c|}{$B_1$}\\
\cline{2-7}N & $L^1$ error&order&  $L^1$ error&order&  $L^1$ error&order\\
\hline
$30^3$	&3.51E-02&	--   & 	1.93E-01&	--	&3.07E-01&	--	\\ \hline
$40^3$	&2.56E-02&	1.10&	1.43E-01&	1.04&	2.28E-01&	1.46\\ \hline
$50^3$	&1.99E-02&	1.12&	1.09E-01&	1.22&	1.79E-01&	2.03\\ \hline
$60^3$	&1.51E-02&	1.51&	8.30E-02&	1.50&	1.37E-01&	2.33
\\ \hline

\hline
\end{tabular}
\end{center}
\end{table}

\subsection{Indirect Approach for Problem 2}
\begin{figure}[htb]
\begin{center}
\includegraphics[width=3in]{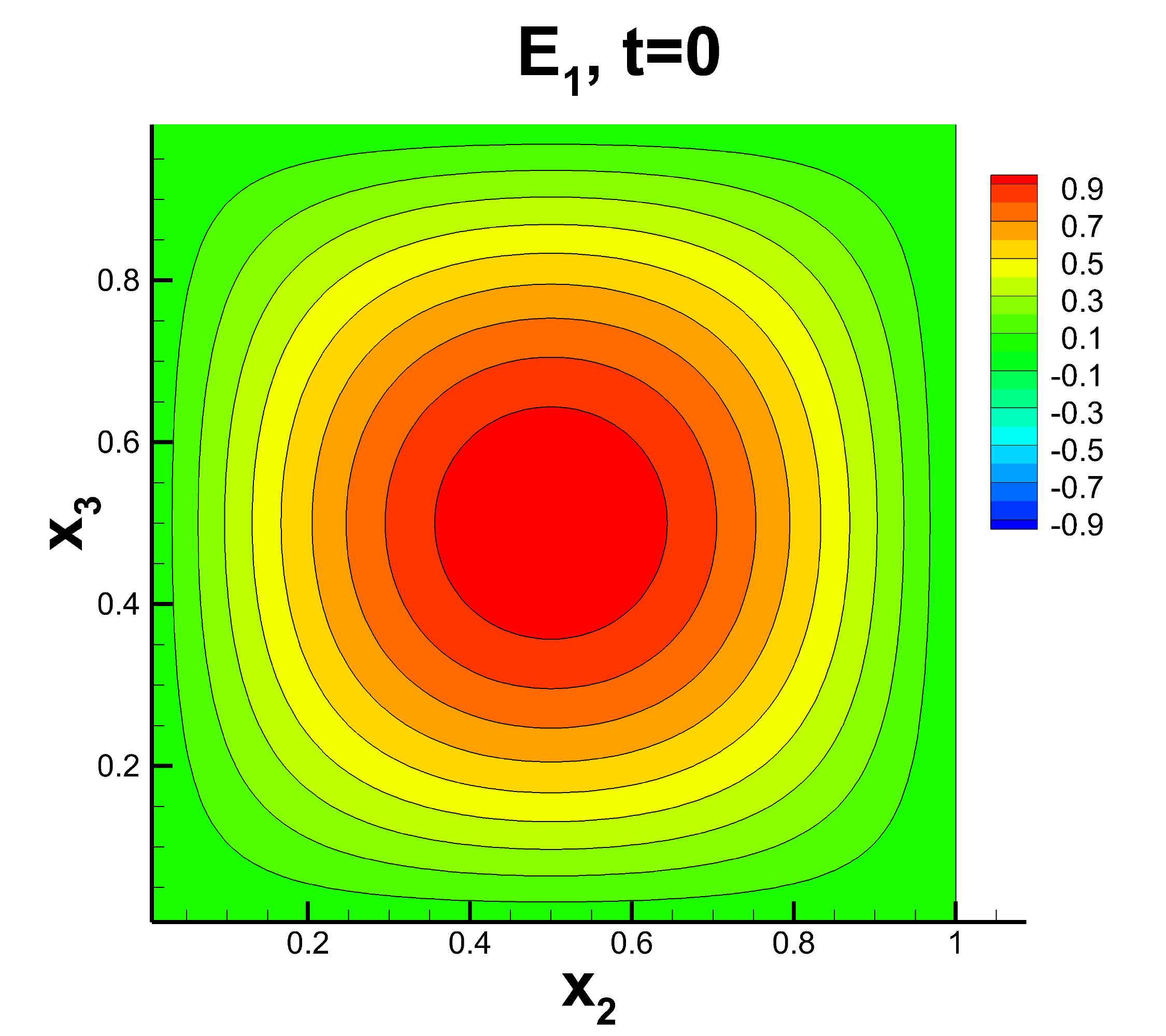}
\includegraphics[width=3in]{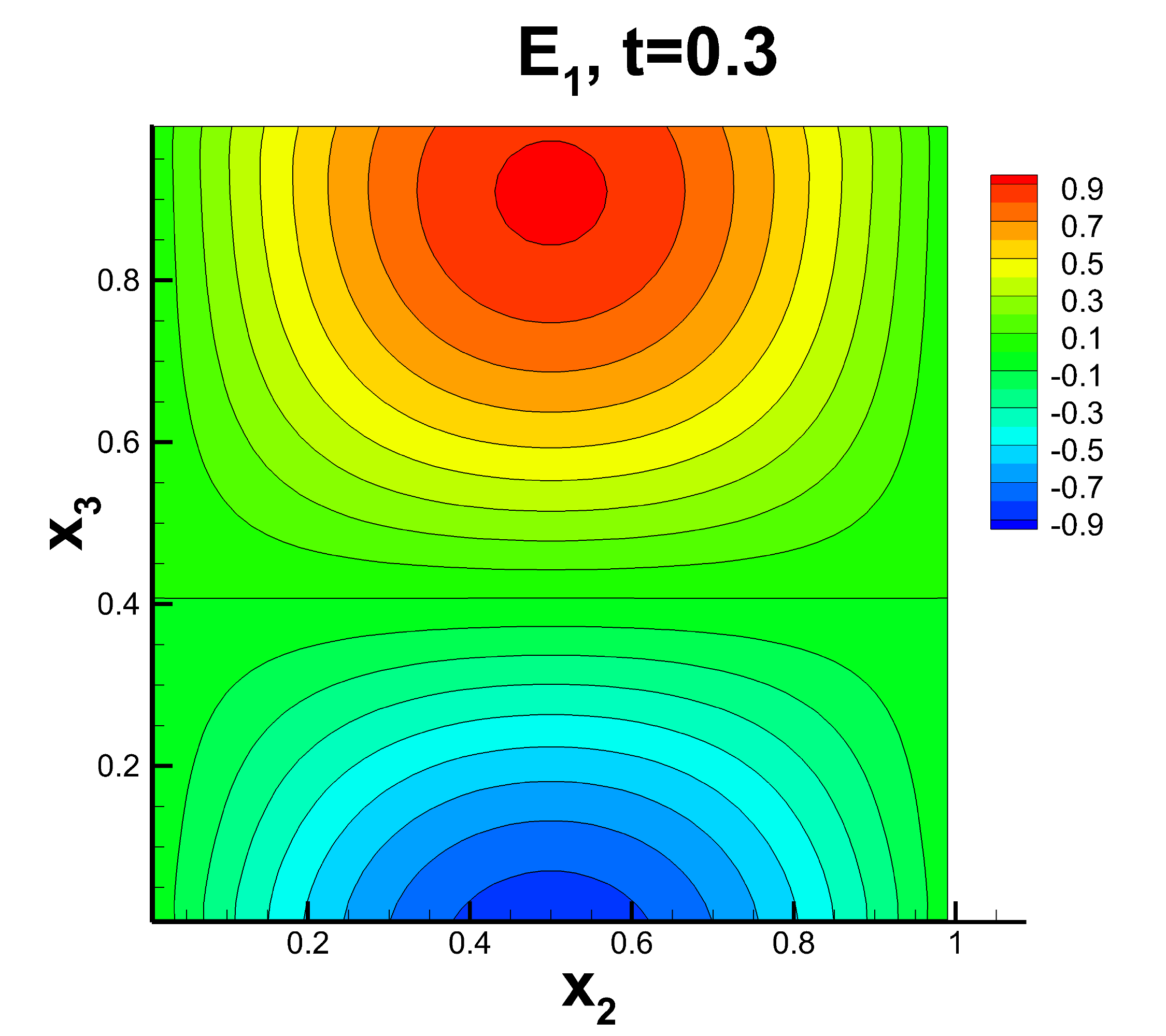}
\includegraphics[width=3in]{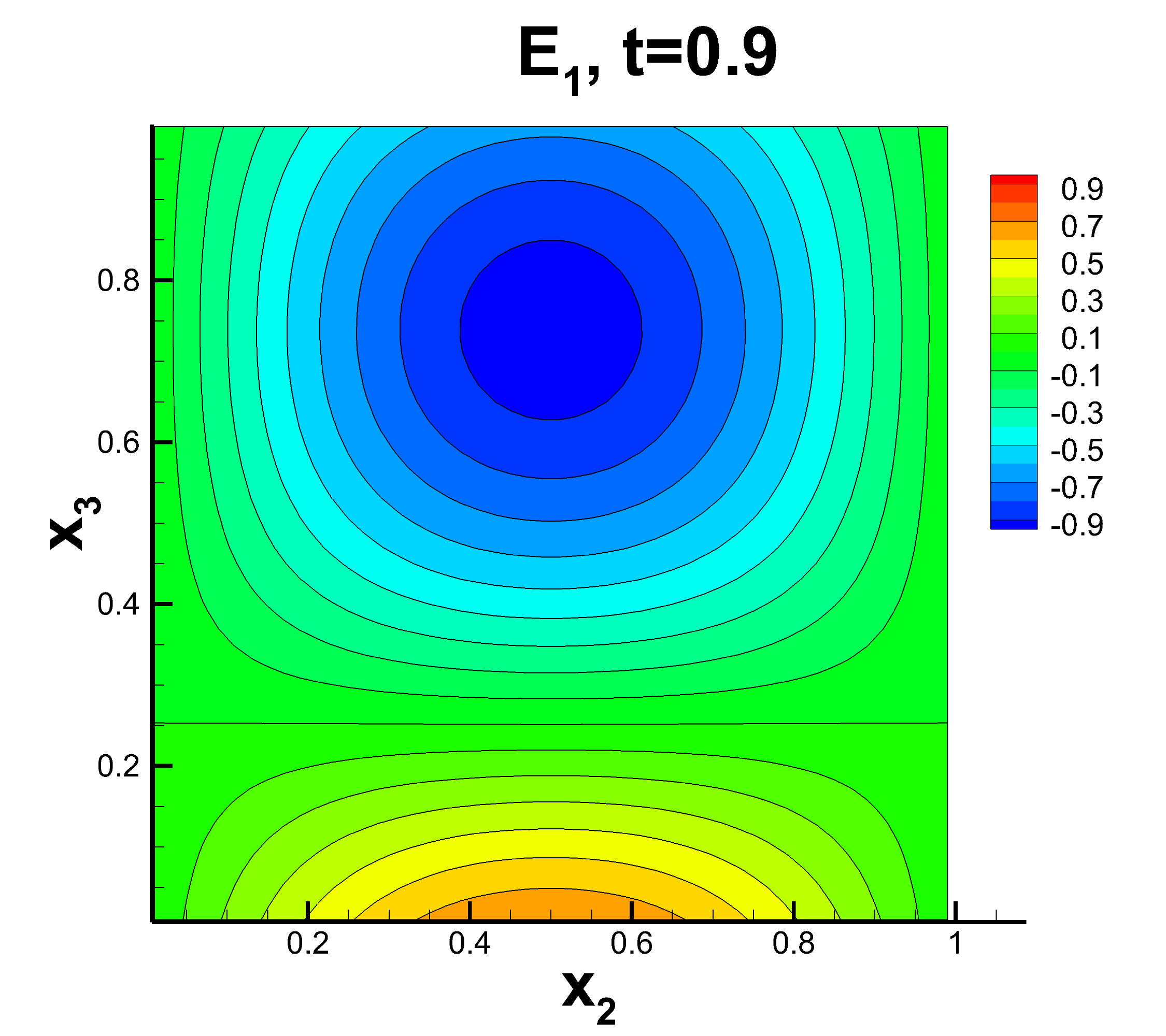}
\includegraphics[width=3in]{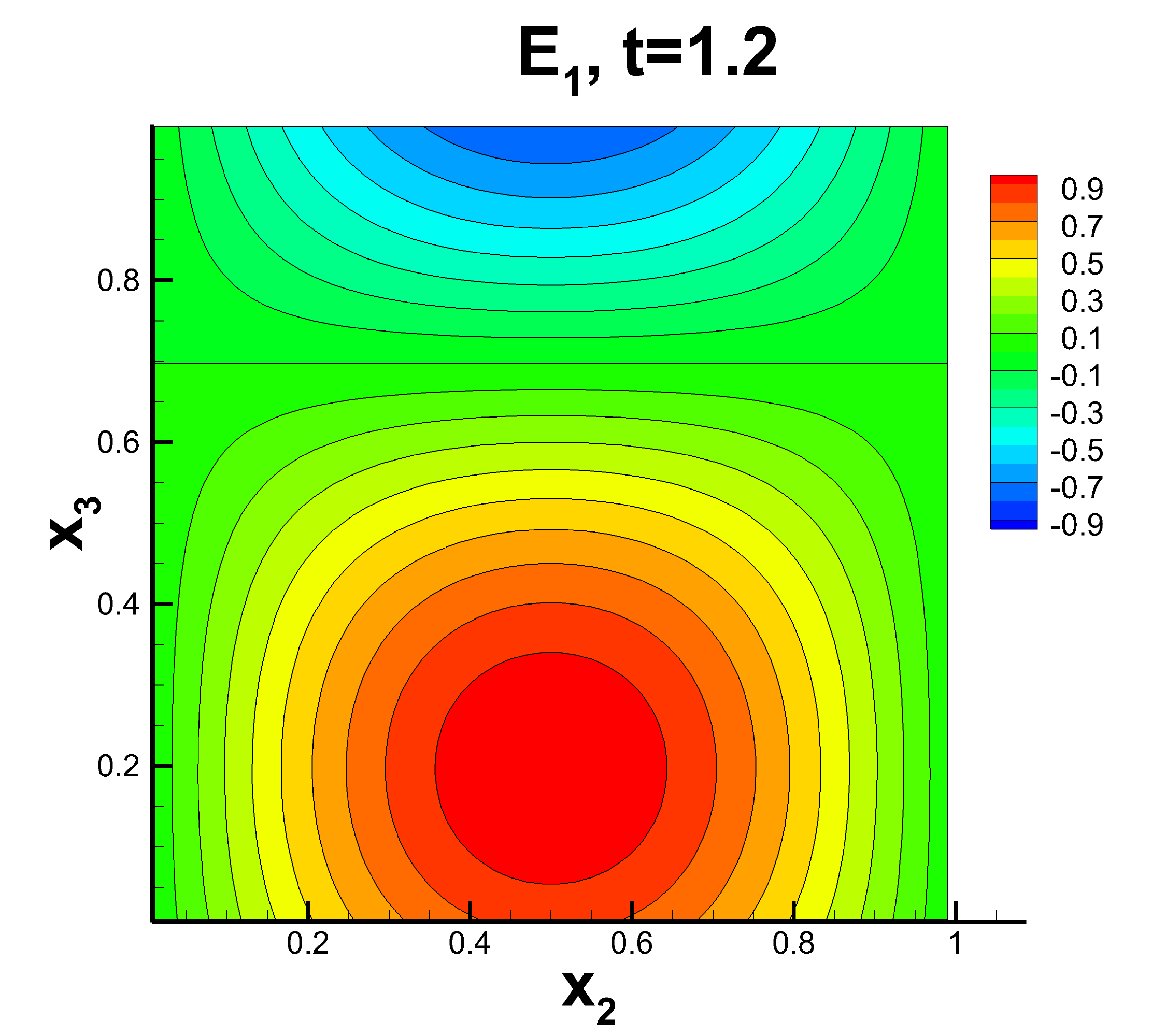}
\end{center}
\caption{
Problem 2.   The contour plots of two-dimensional cuts of $E_1$ at $x_1=0.51$. A TE$_{10}$ mode is propagating in the cubic waveguide. $50^3$ particles. Second order dissipative scheme. Indirect approach. $\text{CFL} =5.0$.
}
\label{fig:sm_eb}
\end{figure}

 Now, we apply the proposed indirect approach for solving problem 2. To save space, we only report the results of the dissipative scheme.
 Table \ref{tab:diss_test1} summarizes the convergence study of the proposed scheme  with a large CFL number 4.9. Second order of convergence is observed for $w_1$, $E_1$ and $B_2$. Several plots of the two-dimensional cuts at $x=0.51$ for the numerical solution $E_1$ are shown in Figure \ref{fig:sm_eb}. Here we use $50^3$ particles in the simulation and let CFL=5. The numerical results are consistent with the exact solution. %At last, we would like to remark again, it is nontrivial to develop decoupled boundary integral formulations for each component of $\bw$, and hence it is more convenient to adopt the second indirect approach if the Silver-M\"{u}ller boundary condition is considered.

\begin{table}[htb]
%\label{tab: tab_sil_mul}
\begin{center}
\caption{Problem 2.      The $L^1$ errors and the corresponding orders of accuracy for $w_1$, $E_1$, and $B_2$. $T=1.5$. Second order dissipative scheme. Indirect approach. CFL=4.9.
\label{tab:diss_test1}}
\bigskip
\begin{tabular}{|c | c c|c c|c c|}
\hline
\cline{2-7} &\multicolumn{2}{c|}{$w_1$} &\multicolumn{2}{c|}{$E_1$} &\multicolumn{2}{c|}{$B_2$}\\
\cline{2-7}N & $L^1$ error&order&  $L^1$ error&order&  $L^1$ error&order\\
\hline
$30^3$	&4.59E-02&	--   & 	2.12E-01&	--	&1.59E-01&	--	\\ \hline
$40^3$	&2.97E-02&	1.50&	1.24E-01&	1.86&	9.46E-02&	1.80\\ \hline
$50^3$	&1.68E-02&	2.56&	7.26E-02&	2.37&	5.43E-02&	2.49
\\ \hline

\hline
\end{tabular}
\end{center}
\end{table}

\subsection{Indirect Approach for Problem 3}
\begin{figure}[htb]
	\begin{center}
		\includegraphics[width=3in]{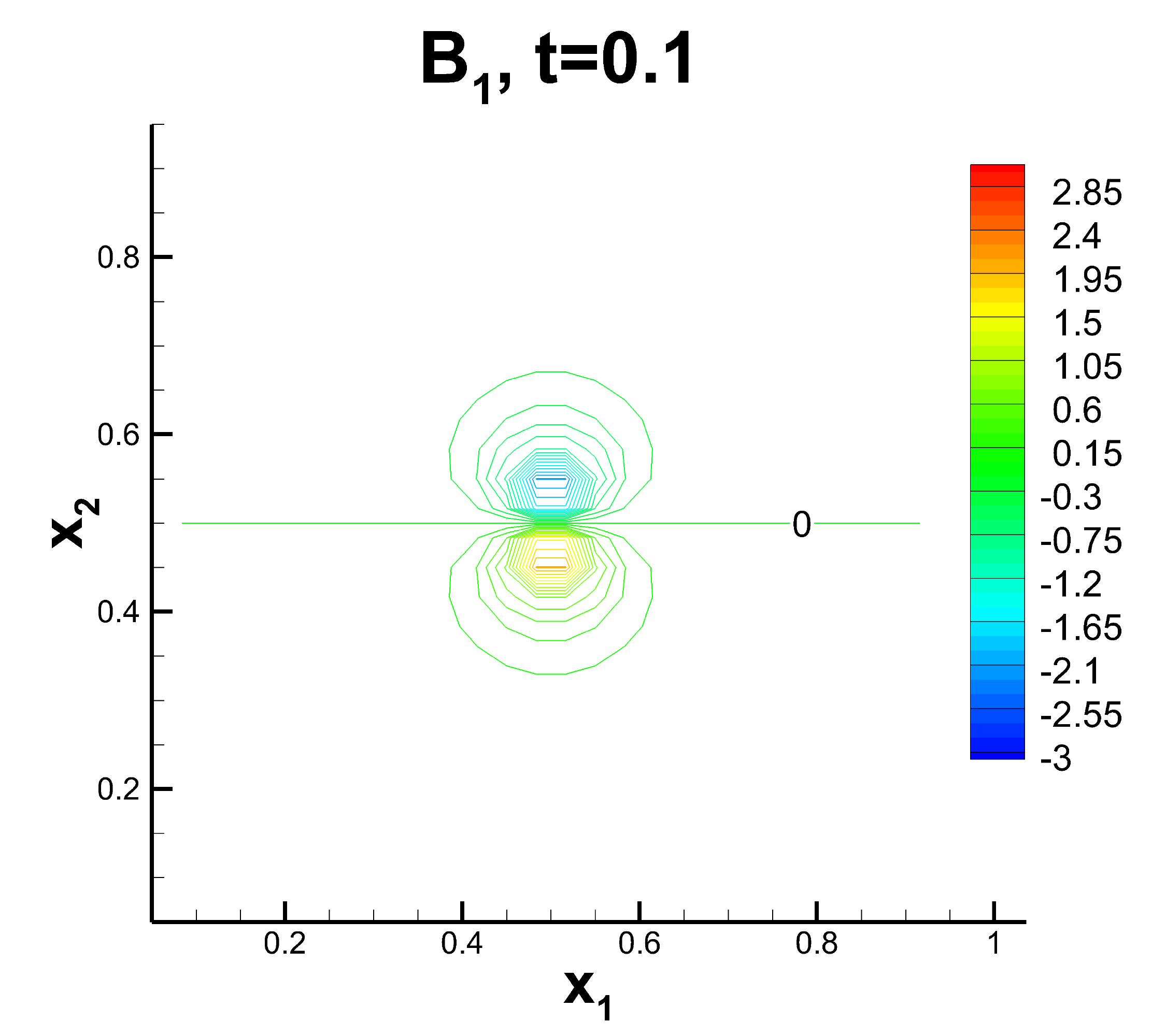}\quad
		\includegraphics[width=3in]{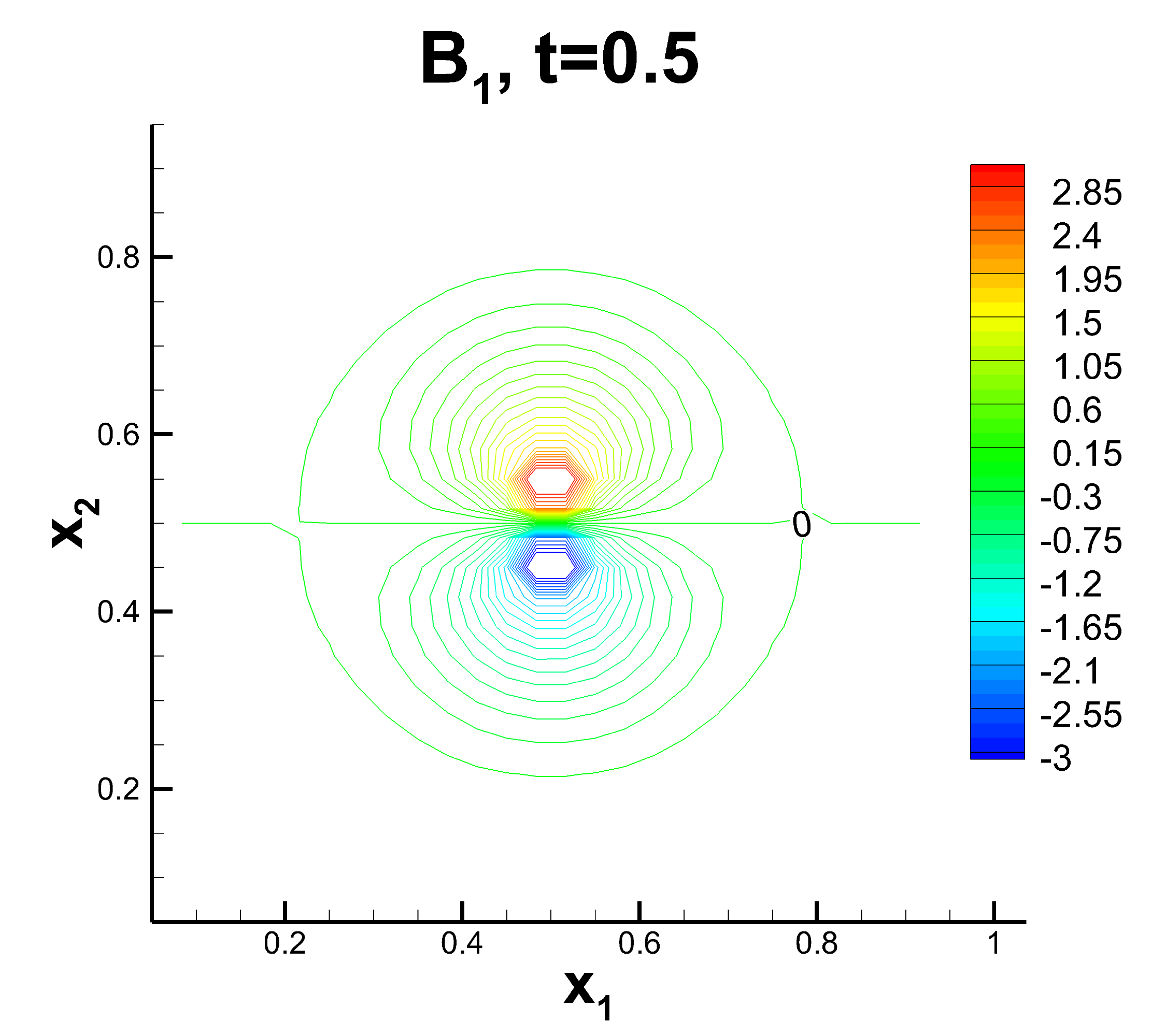}\\
		\includegraphics[width=3in]{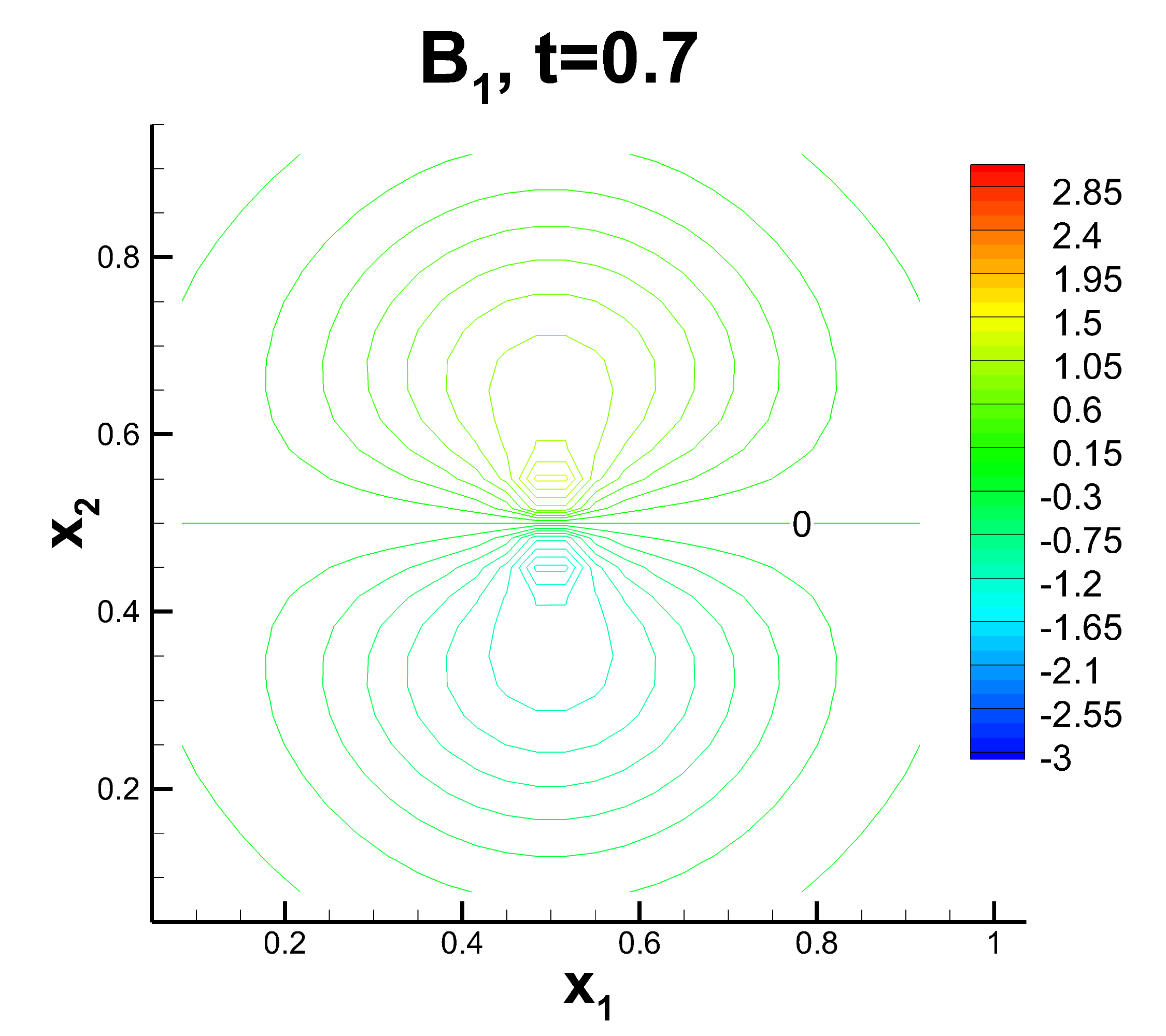}\quad
		\includegraphics[width=3in]{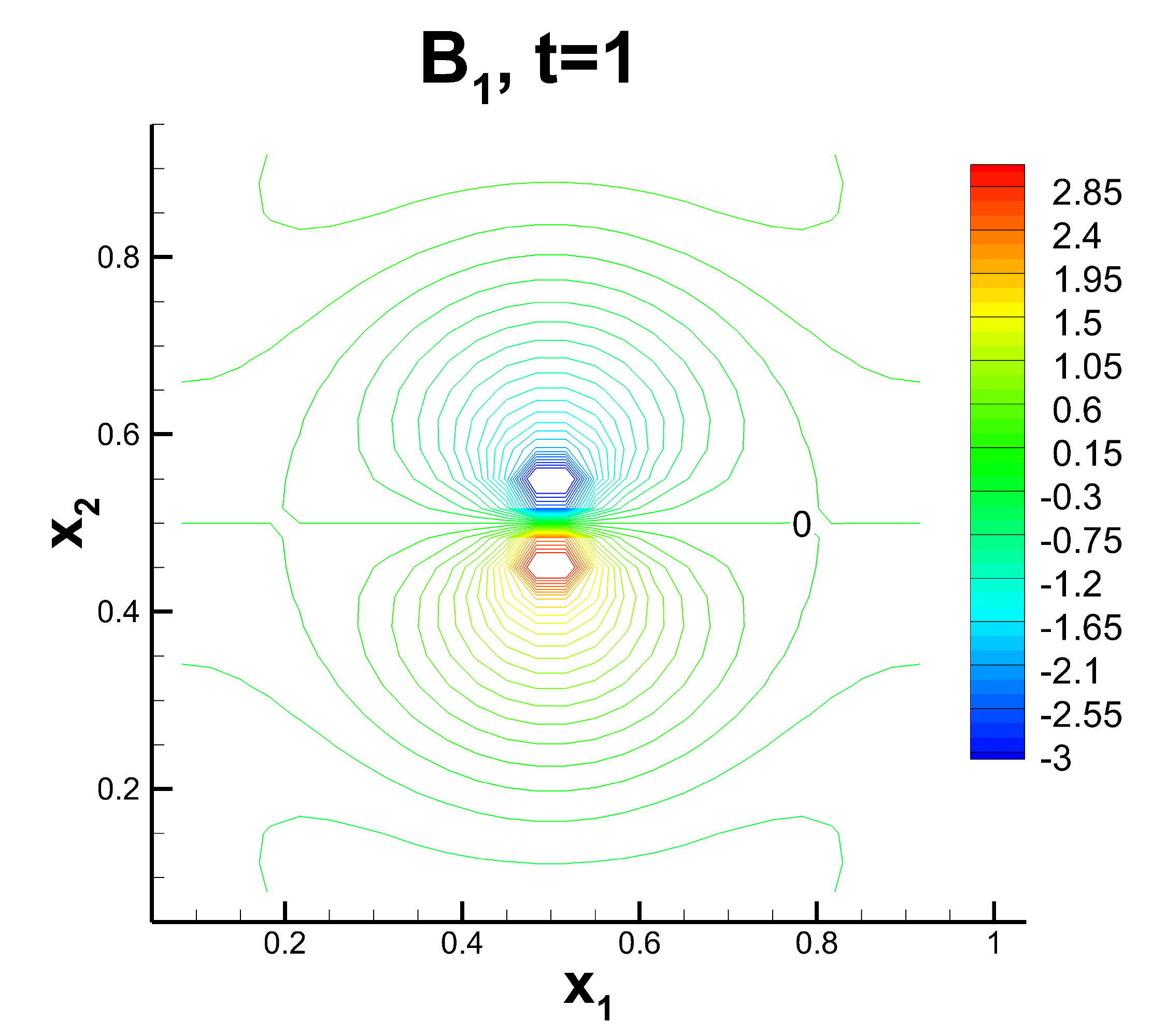}\\
		\includegraphics[width=3in]{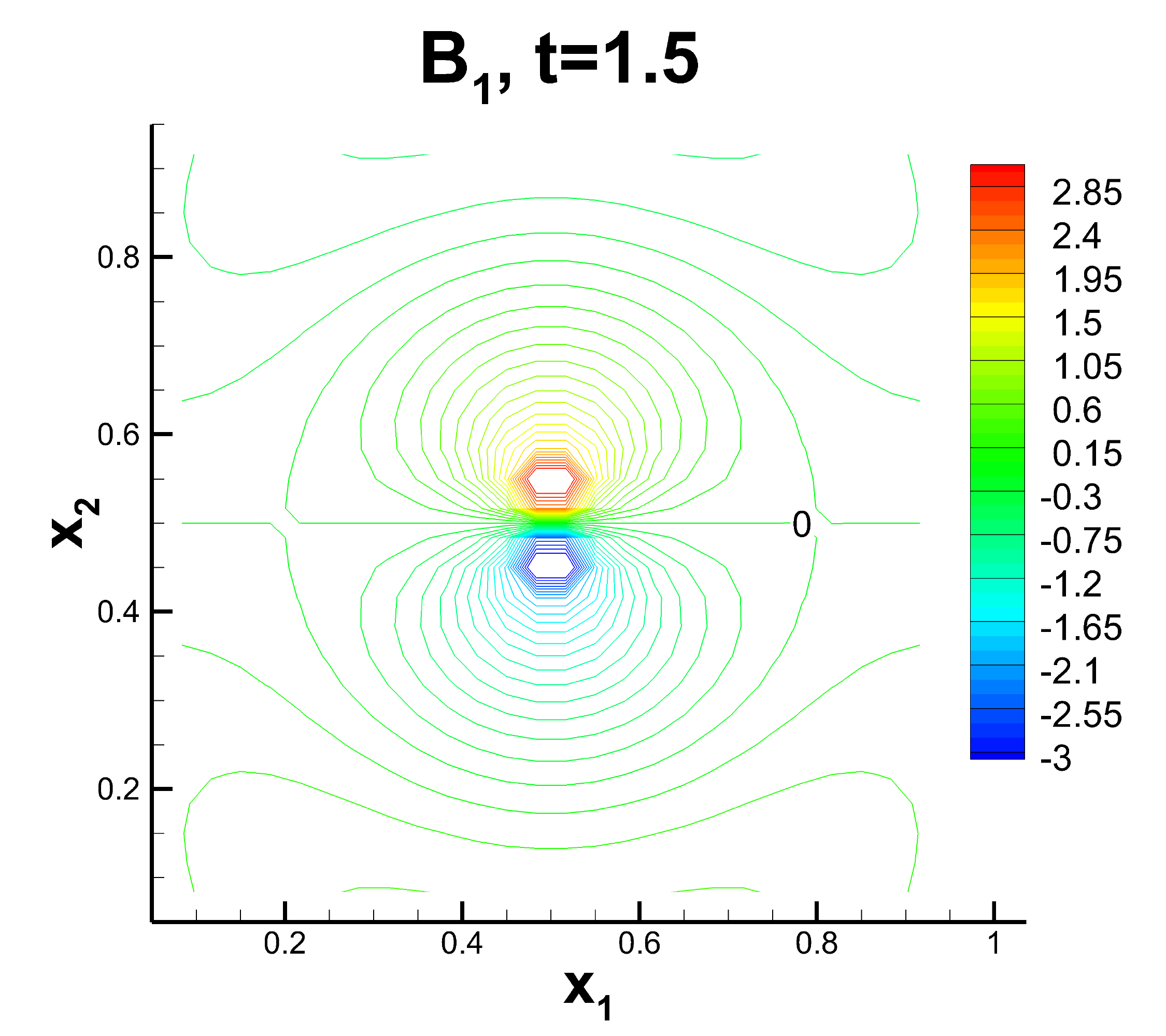}\quad
		\includegraphics[width=3in]{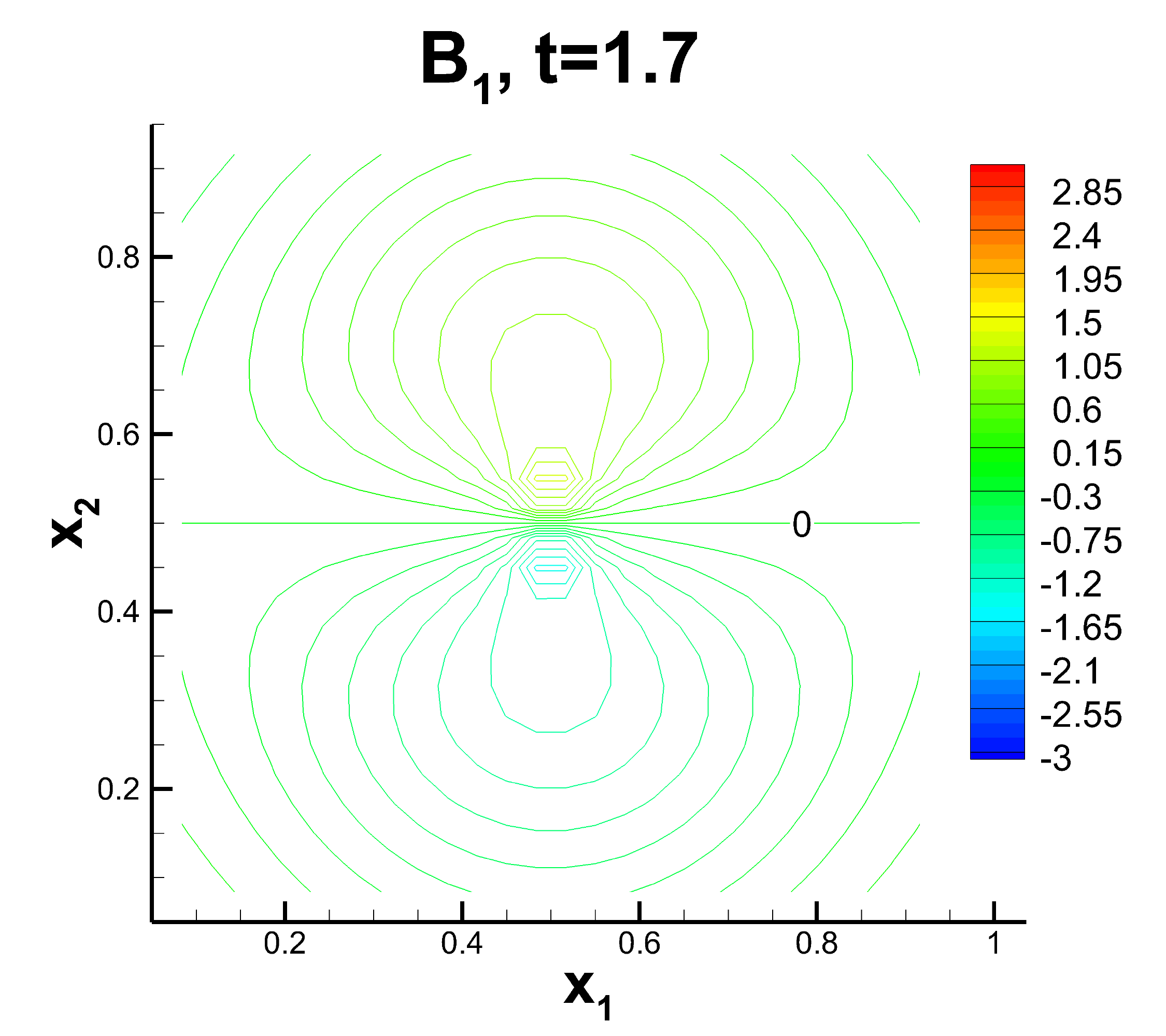}		
	\end{center}
	\caption{
		Problem 3.   The contour plots of two-dimensional cuts of $B_1$ at $x_3=0.2$. $30^3$ particles. Second order dissipative scheme. Indirect approach. $\text{CFL} =3.0$.
	}
	\label{fig:sm_sb}
\end{figure}

Lastly, we present the numerical result for problem 3. Similar to problem 2, we use proposed indirect approach to solve this problem. 
In the simulation, we use $30^3$ particles and let CFL=3. In Figure \ref{fig:sm_sb},  we report the contour plots of two-dimensional cuts of $B_1$ at $x_3=0.2$ at several instances of time. For this problem, the behavior of $\bB$ can be explained by the Biot-Savart law \cite{ciarlet2007continuous}, that is the magnetic field created by the current bar is
$$\bB(\bx)=\frac{\mu_0}{4\pi}\int_C\frac{\bJ\times \mathbf{r}}{|\mathbf{r}|^3}\,dl_{\by},$$
where $\mathbf{r}=\bx-\by$. Therefore, if we denote the orthogonal projection of $\bx$ on the current bar by $\bx'=(1/2,1/2,x_3)$ and let $\mathbf{r}'=\bx-\bx'$, then the Biot-Savart law gives
 $$\mathbf{\bB}(\bx)\propto\frac{\cos(2\pi t)}{|\mathbf{r}'|^2}\left(\begin{array}{c}
 -(x_2-1/2)\\(x_1-1/2)\\0\end{array}
 \right).$$ 
It can be observed from Figure \ref{fig:sm_sb} that at $x_2=1/2$, $B_1$ vanishes horizontally at  $x_2=1/2$. Also note that the Silver-M\"{u}ller absorbing boundary condition is only first order, i.e., only plane waves with normal incidence can be absorbed at the boundary, see \cite{engquist1977absorbing,ciarlet2007continuous}. We can still observe some reflection near the boundary.

\section{Conclusions and future work}
\label{sec:conclusion}

In this paper, we develop AP schemes for Maxwell's equations in the potential form. The methods  are implicit, allow large time steps, and are shown to recover the Darwin limit  at the semi-discrete level when the dimensionless parameter $\epsilon=\bar{v}/c$ goes to 0. By using the  MOL$^T$ framework, we obtain the integral formulation for the vector potential, which are then solved by the treecode algorithm. Although the schemes are only second order accurate in space and time, it is possible to improve the spatial accuracy by using higher order quadrature in the Nystr\"{o}m method framework and  temporal accuracy by the successive convolution technique \cite{causley2014higher}.
Other future directions include   extension of the methods to the scalar and vector potential forms with the Lorentz and Coulomb gauge, as well as incorporation of the schemes in   kinetic plasma simulations.

%For the numerical test, we consider a rectangular waveguide with dimension $[0,1]^3$, in which a TE$_{10}$
% mode propagates in the $z-$direction. The analytical expression of the TE field is given by
%\begin{align*}
%\bE(t)=\left(\begin{array}{c}\sin(\pi y)\sin(\pi z-\omega t)\\
%0\\
%0
% \end{array}\right),\quad
%\bB(t)=\frac{1}{\sqrt{2}}\left(\begin{array}{c}0\\
%\sin(\pi y)\sin(\pi z-\omega t)\\
%\cos(\pi y)\cos(\pi z-\omega t)
% \end{array}\right),\\
%\end{align*}
%where $\omega = \frac{\sqrt{2}\pi}{\varepsilon}$. The exact expression for $\bw$ can be obtained by integrating $\bE(t)$ in time
%\begin{align*}
%\bw(t)=\left(\begin{array}{c}\frac1\omega\sin(\pi y)(\sin(\pi z-\omega t)-\cos(\pi z))\\
%0\\
%0
% \end{array}\right).
%\end{align*}
%We prescribe the perfectly conducting boundary conditions on the four side faces ($x=0$, $x=1$, $y=0$, and $y=1$), while at the bottom ($z=0$) and at the top($z=1$),  the Silver-Muller boundary conditions are imposed, i.e.,
%\begin{align*}
%&\left(\bw_t +\frac{1}{\varepsilon}\nabla\times\bw\times\bn\right)\times\bn = \bg\times\bn,\quad\text{with}\\
%&\bg =\left\{\begin{array}{c} \displaystyle (-\frac{\sqrt{2}+2}{2}\sin(\pi y)\sin(\omega t),\,0,\,\,0)^T\quad\text{on}\quad z=0\\[2mm]
%\displaystyle  (-\frac{\sqrt{2}-2}{2}\sin(\pi y)\sin(\omega t),\,0,\,\,0)^T\quad\text{on}\quad z=1\end{array}\right.
%\end{align*}
%%
%For the numerical simulation, we initialize the potential $\bw$ according to the analytical form at $t=0$.

\bibliographystyle{abbrv}
\bibliography{ref_cheng}

\end{document}